\documentclass{article}
\usepackage{amssymb, latexsym}
\usepackage{fullpage}
\usepackage{amsmath}
\newtheorem{theorem}{Theorem}[section]
\newtheorem{lemma}[theorem]{Lemma}
\newtheorem{lem}[theorem]{Lemma}
\newtheorem{definition}[theorem]{Definition}
\newtheorem{proposition}[theorem]{Proposition}
\newtheorem{prp}[theorem]{Proposition}

\newtheorem{remark}[theorem]{Remark}
\newtheorem{rmk}[theorem]{Remark}
\newtheorem{dfn}[theorem]{Definition}
\newtheorem{deff}[theorem]{Definition}

\newtheorem{cor}[theorem]{Corollary}
\newtheorem{thm}[theorem]{Theorem}

\newtheorem{prop}[theorem]{Proposition}

\newtheorem{rema}[theorem]{Remark}

\newcommand{\comment}[1]{}  %%%to comment out blocks of text				

\newcommand{\step}[3]{
\frac{#1}{#2}\; #3
}

\newcommand{\inn}{\,{\in}\,}

\newcommand{\bI}{{\mathsf{Id}}}

\newcommand{\lando}{\;\wedge\;}
\newcommand{\czf}{\mbox{$\mathbf{CZF}$}}

\newcommand{\fin}{\,{\in}\,}

\newcommand{\ACo}{{\mathbf{AC}}_{\omega}}
\newcommand{\voll}[2]{{\mathbf{mv}}(\phantom{}^{#1}{#2})}

\newcommand{\WW}{{\mathbf W}}
\newcommand{\CST}{{\mathbf{CST}}}

\newcommand{\ist}{=}

\newcommand{\sREA}{\mathbf{sREA}}

\newcommand{\PAs}[1]
{#1^+}

\newcommand{\bpn}{\mathbf{p_0}}
\newcommand{\bpo}{\mathbf{p_1}}

\newcommand{\bes}{\begin{eqnarray*}}
\newcommand{\ees}{\end{eqnarray*}}

\newcommand{\bp}{{\mathbf p}}

\newcommand{\bi}{{\mathbf{i}}}
\newcommand{\bk}{{\mathbf k}}

\newcommand{\PAx}{{\mathbf{PAx}}}
\newcommand{\app}{\mathsf{App}}
\newcommand{\bd}{\mbox{$\mathbf d$}}
\newcommand{\bs}{\mbox{$\mathbf s$}}
\newcommand{\bpw}{\mbox{$\mathbf p_{\mathbf 1}$}}
\newcommand{\dfi}{\;\,:=\;\;}

\newcommand{\I}{\mathbf{I}}
 \newcommand{\PSAC}{\ensuremath{\mathbf{\Pi\Sigma\!-\!AC}}}
\newcommand{\PSIAC}{\ensuremath{\mathbf{\Pi\Sigma\I\!-\!AC}}}
\newcommand{\PSWAC}{\ensuremath{\mathbf{\Pi\Sigma{\mathbf W}\!-\!AC}}}
\newcommand{\PSWIAC}{\ensuremath{\mathbf{\Pi\Sigma{\mathbf W}\I\!-\!AC}}}

\newcommand{\LPO}{{\mathbf{LPO}}}

\def\provx#1#2#3#4{
\setbox1=\hbox{\kern1.5pt$\scriptstyle#3$}
\def\zeichen{#2}
\ifx\zeichen\empty\setbox0=\hbox to .75em{}\else\setbox0=\hbox
{\kern1.5pt$\scriptstyle#2$}\fi
\dimen1=\dp0 \ifdim \dimen1=0pt
\advance \dimen1 by 1.5ex \else \advance \dimen1 by 1.2ex
\fi\dimen3=2ex\dimen4=.5ex\ifdim \wd0<\wd1 \dimen2=\wd1 \else \dimen2=\wd0
\fi\hbox{$#1\hskip 5pt minus5pt\vrule height\dimen3
depth\dimen4\raise\dimen1\copy0\hskip-1\wd0 \lower\ht1
\copy1\hskip-1\wd1\vrule width\dimen2 height.7ex depth-.6ex\hskip3pt
minus1.5pt#4\hskip2pt plus2pt minus2pt$}}

\def\prov#1#2#3{
\setbox1=\hbox{\kern1.5pt$\scriptstyle#2$}
\def\zeichen{#1}
\ifx\zeichen\empty\setbox0=\hbox to .75em{}\else\setbox0=\hbox
{\kern1.5pt$\scriptstyle#1$}\fi
\dimen1=\dp0
\ifdim \dimen1=0pt
\advance \dimen1 by 1.5ex \else \advance \dimen1 by 1.2ex
\fi\dimen3=2ex\dimen4=.5ex\ifdim \wd0<\wd1 \dimen2=\wd1 \else \dimen2=\wd0
\fi\hbox{\hskip0pt plus 4pt
$\vrule height\dimen3
depth\dimen4\raise\dimen1\copy0\hskip-1\wd0
\lower\ht1\copy1\hskip-1\wd1\vrule width\dimen2 height.7ex depth-.6ex
\hskip3pt minus1.5pt#3\hskip2pt plus2pt minus2pt$}}

\def\prv#1#2{
\setbox1=\hbox{\kern1.5pt$\scriptstyle#2$}
\ifx\zeichen\empty\setbox0=\hbox to .75em{}\else\setbox0=\hbox
{\kern1.5pt$\scriptstyle#1$}\fi
\dimen1=\dp0 \ifdim \dimen1=0pt
\advance \dimen1 by 1.5ex \else \advance \dimen1 by 1.2ex
\fi\dimen3=2ex\dimen4=.5ex\ifdim \wd0<\wd1 \dimen2=\wd1 \else \dimen2=\wd0
\fi\hbox{\hskip.5em$\vrule height\dimen3
depth\dimen4\raise\dimen1\copy0\hskip-1\wd0
\lower\ht1\copy1\hskip-1\wd1\vrule width\dimen2 height.7ex depth-.6ex
\hskip3pt minus1.5pt$}}

\mathchardef\str='1066
\def\negprov#1#2#3{
\setbox1=\hbox{\kern1.5pt$\scriptstyle#2$}
\setbox4=\hbox{$\str$}
\def\zeichen{#1}
\ifx\zeichen\empty\setbox0=\hbox to 1em{}\else\setbox0=\hbox
{\kern1.5pt$\scriptstyle#1$}\fi
\dimen1=\dp0
\ifdim \dimen1=0pt
\advance \dimen1 by 1.5ex \else \advance \dimen1 by 1.2ex
\fi\dimen3=2ex\dimen4=.5ex\ifdim \wd0<\wd1 \dimen2=\wd1 \else \dimen2=\wd0
\fi\hbox{\hskip.5em$\kern-1.9pt\raise1pt\copy4\kern-\wd4\kern1.9pt\vrule height\dimen3
depth\dimen4\raise\dimen1\copy0\hskip-1\wd0
\lower\ht1\copy1\hskip-1\wd1\vrule width\dimen2 height.7ex depth-.6ex
\hskip3pt minus1.5pt#3\hskip2pt plus2pt minus2pt$}}

\def\goed#1{\setbox5=\hbox{$#1$}\dimen1=.25em \dimen2=\dimen1 \advance \dimen2
by -1pt\hbox{\raise.65\ht5 \hbox{\vrule height.5\ht5 depth0pt width.4pt\vrule
height.5\ht5 width\dimen1 depth-.48\ht5}\kern-\dimen2\copy5\kern-\dimen2
\raise.65\ht5 \hbox{\vrule height .5\ht5 width\dimen1 depth-.48\ht5\vrule
height.5\ht5 depth 0pt width.4pt}\hskip4pt plus2pt minus2pt}}

\def\mod#1#2{
\def\zeichen{#1}
\hbox{\hskip 2pt plus3pt minus 2pt\vrule width.5pt height2ex depth.5ex
\vbox{\ifx\zeichen\empty\hbox to .75em{}\else
\hbox{\kern1.5pt $\scriptstyle#1$}\fi
\kern2pt
\hrule
\kern1.7pt
\hrule\kern1.7pt}
\hskip3pt minus 2pt$#2$}\hskip2pt
plus3pt minus2pt}

\def\notmod#1#2{\hbox{\hskip 2pt plus 3pt minus 3pt\vrule width.5pt
height2ex depth.5ex
\vbox{\hbox{\kern1.5pt $\scriptstyle#1$}\kern3pt
\setbox0=\hbox{\kern2pt$\scriptstyle/$}
\hrule
\kern-1.7pt
\copy0
\kern-\ht0
\kern 1.7pt
\hrule\kern1.7pt}n
\hskip3pt minus 2pt$#2$}\hskip2pt
plus3pt minus2pt}
\def\sq{\hbox{\rlap{$\sqcap$}$\sqcup$}}
\def\qed{\ifmmode\sq\else{\unskip\nobreak\hfil\penalty50\hskip1em\null
\nobreak\hfil\sq\parfillskip=0pt\finalhyphendemerits=0\endgraf}\fi\medskip}

\def\lleq{\hbox{\hskip3pt minus3pt\kern1pt\lower4pt
\vbox{\hbox{$\scriptstyle\ll$}
\kern-7pt\hbox{\kern1pt$\scriptstyle=$}}\hskip3pt minus 3pt}}

\mathchardef\res='1152
\mathchardef\qin='1062
\mathchardef\qprec='1036
\mathchardef\qless='474
\mathchardef\dpkt='72

\newcommand{\HoTT}{\mathrm{HoTT}}
\newcommand{\TO}{\mathbf{T}_0}
\newcommand{\UA}{\mathbf{UA}}
\newcommand{\tnm}{\TO+\mathbf{MID}}
\newcommand{\ext}{\mbox{\tiny ext}}
\newcommand{\PAxx}{\mathbf{PAx}}
\newcommand{\Inac}{\mathbf{Inacc}}
\newcommand{\FTAC}{\mathbf{FT}\mbox{-}\mathbf{AC}}

\newcommand{\Idd}{\mbox{\sf Id}}

\newcommand{\dekordfuenf}[5]{\begin{array}{l}{#1}\\ {#2}\\ {#3}\\ {#4}\\
{#5}\end{array}}
\newcommand{\dekorddrei}[3]{\begin{array}{l}{#1}\\{#2}\\{#3}\end{array}}

\newcommand{\idpeelJ}{\mbox{\sf J}}
\newcommand{\intt}{\,:\,}
\newcommand{\refl}{\mathsf{refl}}
\newcommand{\CZFe}{\CZF_{Exp}}
\newcommand{\BYst}{\ensuremath{\mathbf{Y^\ast}}}

\newcommand{\lt}{{\mathfrak L}}
\newcommand{\bsn}{\mbox{${\mathbf s}_{\mathbf N}$}}
\newcommand{\sequ}[3]{#1_#2,\ldots,#1_#3}
\newcommand{\brak}[1]{\ulcorner#1\urcorner}
\newcommand{\ltx}{\prec_X}

\newcommand{\Iari}{$\mbox{\bf IA\,\!RI}$ }
\newcommand{\Iario}{$\mbox{\bf IA\,\!RI}$}
\newcommand{\WPU}{\WP_{U}}
\newcommand{\wof}{\mbox{WF}}
\newcommand{\progU}{\prog_{U}}
\newcommand{\prog}{\mathsf{Prog}}

\newcommand{\Proof}{\bx{Proof.}}

\newcommand{\QED}{\mbox{}\hfill$\Box$}	%% requires `latexsym' package
\newcommand{\beq}{\vspace{-2ex}\begin{equation}}
\newcommand{\eeq}{\end{equation}}
\newcommand{\bde}{\begin{definition}}
\newcommand{\ede}{\end{definition}}
\newcommand{\ble}{\begin{lemma}}
\newcommand{\ele}{\end{lemma}}
\newcommand{\bpr}{\begin{proposition}}
\newcommand{\epr}{\end{proposition}}
\newcommand{\bthe}{\begin{Th}}
\newcommand{\ethe}{\end{Th}}
\newcommand{\bco}{\begin{Cor}}
\newcommand{\eco}{\end{Cor}}
\newcommand{\bcl}{\begin{Cla}}
\newcommand{\ecl}{\end{Cla}}
\newcommand{\becon}{\begin{Con}}
\newcommand{\econ}{\end{Con}}
\newcommand{\beoq}{\begin{OQ}}
\newcommand{\eoq}{\end{OQ}}

\newcommand{\ix}[1]{\textit{#1}}

\newcommand{\bx}[1]{\textbf{#1}}

\newcommand{\BY}{\ensuremath{\mathbf{Y}}}

\newcommand{\Ext}{\mathbf{Ext}} %{\ensuremath{\mathrm{Ext}}}

\newcommand{\ZF}{\ensuremath{\mathbf{ZF}}}
\newcommand{\ZFC}{\ensuremath{\mathbf{ZFC}}}
\newcommand{\KP}{\ensuremath{\mathbf{KP}}}

\newcommand{\AC}{\ensuremath{\mathbf{AC}}}

\newcommand{\DC}{\ensuremath{\mathbf{DC}}}
\newcommand{\RDC}{\ensuremath{\mathbf{RDC}}}

\newcommand{\KPi}{\ensuremath{\mathbf{KPi}}}

\newcommand{\CZF}{\ensuremath{\mathbf{CZF}}}

\newcommand{\REA}{\ensuremath{\mathbf{REA}}}

\newcommand{\MLTT}{\ensuremath{\mathbf{MLTT}}}

\newcommand{\MID}{\ensuremath{\mathbf{MID}}}
\newcommand{\UMID}{\ensuremath{\mathbf{UMID}}}

\newcommand{\HA}{\ensuremath{\mathbf{HA}}}

\newcommand{\POT}{\ensuremath{\Pi^1_2}}

\newcommand{\om}{\ensuremath{\omega}}

\newcommand{\Si}{\ensuremath{\Sigma}}

\newcommand{\bj}{\mathbf{j}}

\newcommand{\N}{N}
\newcommand{\A}{\forall}

\newcommand{\lin}{\hspace{-.2em}\in\hspace{-.2em}}

\newcommand{\ra}{\ensuremath{\rightarrow}}

\newcommand{\two}{\uppercase{ii}}

\newcommand{\vin}{\mathbin{\varepsilon}}

\newcommand{\ba}{\begin{array}}
\newcommand{\ea}{\end{array}}

\newcommand{\pair}[2]{{\langle{#1},{#2}\rangle}}

\newcommand{\set}[2]{\{{#1}\,|\:{#2}\}}

\newsavebox{\sbctext}
\newenvironment{ctext}[1]%
	{\begin{lrbox}{\sbctext}\begin{minipage}{#1\unitlength}\centering}%
	{\end{minipage}\end{lrbox}\makebox(0,0){\usebox{\sbctext}}}
\newcommand{\bc}{\begin{ctext}{400}}
\newcommand{\ec}{\end{ctext}}
\newtheorem{Def}{Definition}[section]
\newtheorem{Th}{Theorem}
\newtheorem{Cor}[Def]{Corollary}

\newcommand{\Clop}{\ensuremath{\mathbf{Clop}}}
\newcommand{\Mon}{\ensuremath{\mathbf{Mon}}}
\newcommand{\Lfp}{\ensuremath{\mathbf{Lfp}}}
\newcommand{\OP}{\ensuremath{\mathbf{OP}}}

\newcommand{\topp}{{\mathbf{T}}^{\OP}_{<\omega}}
\newcommand{\Topp}{{\mathbf{T}}^{\OP}_{<\varepsilon_0}}

\newcommand{\bN}{{\mathbb N}}
\newcommand{\beqs}{\begin{eqnarray*}}
\newcommand{\eeqs}{\end{eqnarray*}}
\newcommand{\beqq}{\begin{eqnarray}}
\newcommand{\eeqq}{\end{eqnarray}}

\newcommand{\oder}{\;\vee\;}

\newcommand{\paar}[1]{\langle#1\rangle}

\newcommand{\lra}{\leftrightarrow}

\newcommand{\WP}{{\mathbf{WP}}}

\newcommand{\ini}{{\in}}

\newcommand{\then}{\rightarrow}
\newcommand{\prf}{{\bf Proof\/}: }

\newcommand{\tn}{\mbox{${\mathbf T}_{\mathbf 0}$}}

\newcommand{\tnurc}{\mbox{${\mathbf T}_{\mathbf 0}^c\!\restriction+\UMID_{\bN}$}}
\newcommand{\tnurindc}{\mbox{${\mathbf T}_{\mathbf 0}^c\!\restriction+\indd+\UMID_{\bN}$}}
\newcommand{\tnuri}{\mbox{${\mathbf T}_{\mathbf 0}^i\!\restriction+\UMID_{\bN}$}}
\newcommand{\tnurindi}{\mbox{${\mathbf T}_{\mathbf 0}^i\!\restriction+\indd+\UMID_{\bN}$}}

\newcommand{\indd}{{\mathbf{IND}}_{\mathbb N}}
\newcommand{\clfp}{{\mathbf{lfp}}}
\newcommand{\pizwei}{({\mathbf{\Pi}}^1_2-{\mathbf{CA}})}
\newcommand{\pizweir}{({\mathbf{\Pi}}^1_2-{\mathbf{CA}})_0}

\newcommand{\dg}{=_{\!\mbox{\tiny ext}}}
\newcommand{\elee}{\varepsilon}
\newcommand{\extsubseteq}{\subseteq}
\newcommand{\V}{V}

\begin{document}
%\title{The Role of Constructive and Semi-Constructive Systems in Feferman's Work}
\title{Proof Theory of Constructive Systems: Inductive Types and Univalence}

\author{Michael Rathjen \\ Department of Pure Mathematics\\
 University of Leeds\\ Leeds LS2 9JT, England \\ rathjen@amsta.leeds.ac.uk}

\date{}

\maketitle

\begin{abstract}   In Feferman's work, explicit mathematics and theories of generalized inductive definitions  play a central role. One objective of this article is to describe the connections with Martin-L\"of type theory and constructive Zermelo-Fraenkel set theory. Proof theory has contributed to a deeper grasp of the relationship between different frameworks for constructive mathematics. Some of the reductions are known only through ordinal-theoretic characterizations. The paper also addresses the strength of Voevodsky's univalence axiom.
%which will be pointed out in the article.

A further goal is to investigate the strength of intuitionistic theories of generalized inductive definitions in the framework of intuitionistic explicit mathematics that lie beyond the reach
of Martin-L\"of type theory.   \\
{\it Key words: Explicit mathematics, constructive Zermelo-Fraenkel set theory,  Martin-L\"of type theory, univalence axiom, proof-theoretic strength} \\
MSC 03F30 03F50 03C62
 \end{abstract}

\section{Introduction}
Intuitionistic systems of inductive definitions have figured prominently in Solomon Feferman's program of reducing classical subsystems of analysis and
theories of iterated inductive definitions to constructive                                                                                                                                                                                                                                                                                                                                                                                                                                                                                                                                                                                              theories of various kinds.
  In the special case of classical theories of finitely as well as transfinitely iterated inductive definitions, where the iteration occurs along a computable well-ordering,                                      the program was mainly completed by Buchholz, Pohlers, and Sieg more than 30 years ago (see \cite{BFPS,Fef2010}).
  For stronger theories of inductive definitions such as those based on Feferman's intutitionic {\em Explicit Mathematics}\footnote{Feferman introduced the theory of explicit mathematics in \cite{Fef75}. There it was based on intuitionistic logic and
  notated by $\mathrm{T}_0$. The same notation is used e.g. in \cite{BFPS,j83,rohio} but increasingly $\mathrm{T}_0$ came to be identified with its  classical version. As a result, we adopt the notation $\TO^i$ to stress its intuitionistic basis and reserve $\TO$ for the classical theory.}
  ($\TO^i$) some answers have been provided in the last 10 years while some questions are still open.

  The aim of the first part of this paper is to survey the landscape of some prominent constructive theories that emerged in the 1970s. In addition to Feferman's  $\TO^i$,  Myhill's {\em Constructive Set Theory} ($\CST$) and  Martin-L\"of type theory ($\MLTT$) have been proposed with the aim of isolating the principles on which constructive
mathematics is founded, notably the notions of constructive function and set in Bishop's mathematics.

Martin-L\"of type theory with
   infinitely many universes and inductive types ($\mathsf W$-types) has attracted a great deal of attention recently because of a newly found
   connection between type theory and topology, called {\em homotopy type theory} ($\HoTT$), where types are interpreted as spaces, terms as maps and the inhabitants of the iterated identity types on a given type $A$ are viewed as paths, homotopies and higher homotopies of increasing levels, respectively, endowing each type with a weak $\omega$-groupoid structure.

Homotopy type theory, so it appears, has now reached the mathematical mainstream:
 \begin{quote}\label{hybris}{\em  Voevodsky's Univalent Foundations require not just one inaccessible cardinal but an infinite string of cardinals, each inaccessible from its predecessor.} (M. Harris, {\em Mathematics without apologies}, 2015).\end{quote}

 By Univalent Foundations Harris seems to refer to $\MLTT$ plus Voevodsky's {\em Univalence Axiom} ($\UA$).
 To set the stage for the latter axiom, let us recall a bit of history of extensionality  and universes in type theory. 
  Simple type theory, as formulated by A. Church in 1940 \cite{church}, already provides a natural and elegant  alternative to set theory for representing mathematics in a formal way. The stratification of mathematical objects into the types of propositions, individuals and functions between two types is indeed quite natural. In this setup, the axiom of extensionality comes in two forms: the stipulation that two logically equivalent propositions are equal and the stipulation that two pointwise equal functions are equal.
  Some restrictions of expressiveness encountered in simple type theory are overcome by dependent type theory, yet  still  unnatural limitations remain  in that one cannot express the notion of an arbitrary structure in this framework.
   For instance one cannot assign a type to an arbitrary field. Type theory (and other frameworks as well) solve this issue by introducing the notion of a {\em universe type}. Whereas most types come associated with a germane axiom of extensionality inherited from its constituent types following the example of simple type theory, it is by no means clear what kind of extensionality principle should govern universes. A convincing proposal  was missing until the work of V. Voevodsky with its formulation of the extensionality axiom for universes in terms of equivalences. This is the univalence axiom, which generalizes propositional extensionality.
 
 Harris's claim that an infinite sequence of inaccessible cardinals is required to model $\MLTT$ plus Voevodsky's {\em Univalence Axiom} is a pretty strong statement.
 Recent research by Bezem, Huber, and Coquand (see \cite{bch}), though, indicates that $\MLTT+\UA$ has an interpretation in $\MLTT$
 and therefore is proof-theoretically not stronger than $\MLTT$.
 But what is the strength of $\MLTT$?
 As there doesn't seem to exist much common knowledge among type theorists about the strength of various systems and how they relate to the other constructive frameworks as well as classical theories used as a classification hierarchy in reverse mathematics
 and set theory, it seems reasonable to devote a section to mapping out the relationships and gathering current knowledge
 in one place. In this section attention will also be payed to the methods employed in proofs such as interpretations
but  with a particular eye toward the role
 of ordinal analysis therein.

The second part of this paper (Section \ref{secmon}) will be concerned with extensions of explicit mathematics by principles that allow the construction of inductive classifications that lie way beyond $\MLTT$'s reach but still have a constructive flavor.
The basic theory here is intuitionistic explicit mathematics $\TO^i$.  In $\TO^i$ one can
freely talk about monotone operations on classifications and assert the existence of least fixed points of such operators.
 There are two ways in which one can add a principle to $\TO^i$ postulating the existence of least fixed points.
 $\MID$ merely existentially asserts that every monotone operation has a least fixed point whereas $\UMID$ not only postulates the
existence of a least solution, but,
 by adjoining a new functional  constant to the language,
ensures that a fixed point is uniformly
presentable as a function of the monotone operation.

The question of the strength of  %$\tn+\MID$ and $\tn+\UMID$
systems of explicit mathematics with $\MID$ and $\UMID$ was
raised by Feferman in \cite{Fef82}; we quote:
\begin{quote}{\em What is the strength of $\tnm$? [...] I have tried, but did
not succeed, to extend my interpretation of $\tn$ in $\Sigma^1_2-AC + BI$ to
include the statement $\MID$.
The theory $\tnm$ includes all constructive
formulations of iteration of monotone inductive definitions of which I am
aware, while $\tn$ (in its $IG$ axiom) is based squarely on the general
iteration of accessibility inductive definitions. Thus it would be of great
interest for the present subject to settle the relationship between these
theories.} (p. 88)
\end{quote}
As it turned out, the principles $\MID$ and even more $\UMID$ encapsulate considerable strength, when considered on
the basis of classical $\TO$. For instance $\TO+\UMID$ embodies the strength of $\Pi^1_2$-comprehension.
The first (significant) models of $\TO+\MID$ were found by Takahashi \cite{tak}. Research on the precise strength was conducted by Rathjen\cite{Ra96,Ra98,Ra99} and Gla{\ss}, Rathjen, Schl\"uter \cite{grs}. The article \cite{Ra02} provides a survey of the classical case.  Tupailo \cite{Tumid} obtained the first result in the intuitionistic setting.
This and further results will be the topic of section 3.

\section{Some Background on Feferman's $\TO^i$}% and CZF}
The theory of explicit mathematics, here denoted by  $\TO^i$, is a formal framework
that has great expressive power.  It is suitable for representing Bishop-style constructive
mathematics as well as generalized recursion, including direct expression of
structural concepts which admit self-application.
Feferman was led to the development of his {\em explicit mathematics} when trying to understand what Errett Bishop
had achieved in his groundbreaking constructive redevelopment of analysis in \cite{bishop}.
For a detailed account see \cite{Fef75,Fef79}.
The ontology behind the axioms of $\TO^i$ is that the universe of mathematical objects is populated by (a) natural numbers, (b) operations (in general partial) and (c) classifications (akin to Bishop's sets) where operations and classifications are to be understood as given intensionally. Operations can be applied to any object including operations and classifications; they are governed by axioms giving them the structure of a {\em partial combinatory algebra} (also known as {\em applicative structures} or  {\em Sch\"onfinkel algebras}).  There are, for example, operations that act on classifications $X,Y$ to produce their Cartesian product $X\times Y$ and exponential $X^Y$. The formation of classifications is governed by the {\em Join},
{\em Inductive Generation} and {\em Elementary Comprehension Axiom}.

The language of  $\TO^i$,  $\lt(\TO^i)$, has two
sorts of variables. The
free and bound variables $(a,b,c,\ldots$ and $x,y,z\ldots)$
are conceived to range over the whole constructive universe which
comprises {\it operations\/} and {\it classifications\/} among other kinds of
entities; while upper-case versions of
these
$A,B,C,$ ...  and  $X,Y,Z$, ... are
used to represent free and bound classification variables.

{\bf N}  is a classification constant taken to define
the class of natural numbers.  {\bf 0} ,
${\mathbf s}_{\mathbf N}$  and
${\mathbf p}_{\mathbf N}$ are operation constants whose intended
interpretations are the natural number  0  and the successor and
predecessor operations.  Additional
operation constants are  {\bf k}, {\bf s}, {\bf d},
{\bf p}, $\bpo$ and $\bpw$
for the two basic combinators,
definition by cases on  {\bf N}, pairing and the corresponding two
projections.
Additional classification
constants are generated using the axioms and the constants
{\bf j}, {\bf i}  and ${\mathbf c}_{\mathbf n} ({\mathbf {n}} < {\mathbf {\omega}})$
for {\it join}, {\it
induction}  and  {\it comprehension}.

There is no arity associated with the various constants.  The {\it terms} of
$\TO^i$  are just the variables and
constants of the two sorts.  The atomic
formulae of  $\TO^i$  are
built up using the terms and three primitive
relation symbols  =, $\app$  and  $\vin$ as follows.
If $q,r,r_1,r_2$ are terms,
then $q=r$, $\app(q,r_1,r_2)$, and $q\vin r$ (where $r$ has to be a
classification variable or constant) are atomic formulae.
$\app(q,r_1,r_2)$
expresses that the operation $q$ applied to $r_1$ yields the value
$r_2$; $q\vin r$ asserts\footnote{It should be pointed out
that we use the symbol ``$\vin$'' instead of ``$\in$'' deliberately, the latter  being
reserved for the set--theoretic elementhood relation.}
that $q$ is in $r$ or that $q$ is classified
under $r$.

We write $t_1t_2\simeq t_3$ for $\app(t_1,t_2,t_3)$.
%The formulas are obtained by closing under $\land,\then,$ and quantification
%for both sorts.

The set of formulae is then obtained from these using the
propositional connectives and the
two quantifiers of each sort.
%The intended interpretation of
%$\mbox{App}(a,b,c)$  is ``$ a$  applied to  $b$  yields the value  $c$.

In order to facilitate the
formulation of the axioms, the language of $\TO^i$  is
expanded definitionally with the symbol
$\simeq$  and the auxiliary notion of an {\it application
term} is introduced.  The
set of application terms is given by two clauses:
\begin{enumerate}
  \item {all terms of  $\TO^i$  are application terms; and}
  \item if  {$ s$} and  $ t$ are application terms,
then  $( st)$  is an application term.
\end{enumerate}
If $s$ is an application term and $u$ is a bound or free variable we define
$ s\simeq u$  by induction on the buildup of $s$:
$${ s} \simeq { u}\; {\;\mbox{ is }\;}\begin{cases}
                     { s} = { u},& \text{if $ s$ is
a variable or a constant,}\\
\exists x,y [{s}_{ 1}
\simeq x \,\wedge\,{ s}_{ 2} \simeq y
\,\wedge\, \app(x,y,{ u}] &
\text{if $ s$ is an application term $(s_1s_2)$}\end{cases}$$
%}=(${\mathbf s}_{\mathbf 1}$,${\mathbf
%s}_{\mathbf2}$) .}
For  $s$ and  $t$
application terms, we have auxiliary, defined formulae of the form:
\begin{eqnarray*} s \simeq  t& \dfi &\forall y  ( s \simeq y
\leftrightarrow  t \simeq y).\end{eqnarray*}
Some abbreviations are ${ t}_{ 1}\ldots
{ t}_{ n}$ for ((...(${
t}_{1}{ t}_{ 2}$)...)${ t}_{ n}$);
$t\downarrow$ for $\exists y({
t}\simeq y)$ and $\phi ({ t})$
for $\exists y({ t}\simeq y \wedge \phi (y))$.

G\"{o}del numbers for formulae
play a key role in the axioms introducing the classification
constants ${\mathbf c}_{\mathbf n}$.
A formula is said to be {\it elementary} if it contains only free
occurrences of classification variables $A$
(i.e., only as {\it parameters}), and even those free occurrences of $A$ are
restricted: $A$ must occur only to the right of $\vin$ in atomic
formulas. The G\"{o}del number ${\mathbf c}_{\mathbf n}$ above
is  the G\"{o}del number of an elementary formula.
We assume that a standard
G\"odel numbering numbering has been chosen for $\lt(\TO^i)$; if $\phi$ is an
elementary formula and $a,\sequ b1m,\sequ A1n$ %$(n,m>0)$
is a list
of variables which includes all parameters of $\phi$, then
$\{x : \phi(x,\sequ b1n,\sequ A1n)\}$ stands for
$\mathbf c_{\mathbf n}(\sequ b1n,\sequ A1n)$; $\mathbf n$ is the code of the pair of
G\"odel numbers $\langle\brak{\phi}$, $\brak{(a,\sequ b1m,\sequ A1n)}\rangle$
and is
called the
\lq index' of $\phi$ and the list of variables.

Some further conventions are useful.
Systematic notation for $n$-{\it tuples\/} is introduced as follows:
$(t)$ is $t$, $(s,t)$ is $\bp st$, and $(\sequ t1n)$ is defined by
$(( t_1,\ldots,t_{n-1}),t_n)$. Finally, $t'$ is written for the term $\bsn t$, and
$\perp$ is the elementary formula $\mathbf 0\simeq\mathbf 0'$.\\[2ex]
$\TO^i$'s
logic is intuitionistic two-sorted predicate logic with identity.  Its
non-logical axioms are:
\paragraph{I.  Basic Axioms}
\begin{enumerate}
  \item{$\forall X \exists x (X=x)$}
\item $\app(a,b,c_1)\lando\app(a,b,c_2)\,\then\, c_1=c_2$
%  \item{$App({\mathbf a},{\mathbf b},{\mathbf c}_{\mathbf 1})
%\wedge App({\mathbf a},{\mathbf b},{\mathbf c}_{\mathbf 2})
%\rightarrow
%  {\mathbf c}_{\mathbf 1}={\mathbf c}_{\mathbf 2}$}
\end{enumerate}
\paragraph{II.  App Axioms}
\begin{enumerate}

\item $(\bk ab)\downarrow\lando \bk ab\simeq a$,
\item $(\bs ab)\downarrow\lando \bs abc\simeq ac(bc)$,
\item $(\bp a_1a_2)\downarrow\lando(\mathbf p_{\mathbf 1} a)\lando
(\mathbf p_{\mathbf 2} a)\downarrow\lando \mathbf p_{\mathbf i}
(\bp a_1a_2)\simeq a_i$
for $i=0,1$,
\item $(c_1=c_2\lor c_1\ne c_2)\lando (\bd abc_1c_2)\downarrow\lando(c_1=c_2
\then \bd abc_1c_2\simeq a)\lando (c_1\ne c_2\then \bd abc_1c_2\simeq b)$,
\item $a\vin\bN\lando b\vin\bN\then[a'\downarrow\lando \bpn(a')
\simeq a\lando \neg(a'\simeq 0)\lando (a'\simeq b'\then a\simeq b)]$.
\end{enumerate}

%\begin{enumerate}
 % \item{({\bf kab})$\downarrow \wedge$ {\bf kab} $\simeq$ {\bf a}}
  %\item{({\bf sab})$\downarrow \wedge$ {\bf sab} $\simeq$ {\bf ac}({\bf bc})}
%  \item{${\mathbf c}_{\mathbf 1}$ $\in$ {\bf N}
%$\wedge$ ${\mathbf c}_{\mathbf 2}$ $\in$ {\bf N} $\rightarrow$
 % [(${\mathbf {dabc}}_{\mathbf 1}{\mathbf c}_{\mathbf 2}$)
%$\downarrow\wedge$ $(({\mathbf c}_{\mathbf
 % 1}$=${\mathbf {c_2}})\rightarrow {\mathbf d}{\mathbf a}
%{\mathbf b}{\mathbf c}_{\mathbf 1}{\mathbf c}_{\mathbf 2}\simeq$
 % {\bf a})$\wedge$
%  $(({\mathbf c}_{\mathbf 1} \neq {\mathbf c}_{\mathbf 2}$)
%  $\rightarrow {\mathbf d}{\mathbf a}
%{\mathbf b}{\mathbf c}_{\mathbf 1}{\mathbf c}_{\mathbf 2}\simeq {\mathbf b})]$}
%  \item{${\mathbf a}\in {\mathbf N}
%\wedge {\mathbf b} \in {\mathbf N}\rightarrow [{\mathbf a'}\downarrow\wedge{\mathbf
%  p}_{\mathbf N}{\mathbf a'}\downarrow
%\wedge {\mathbf p}_{\mathbf N}({\mathbf a'})\simeq{\mathbf a} \wedge \neg({\mathbf
%  a'}={\mathbf 0}) \wedge ({\mathbf a'}
%\simeq {\mathbf b'}\rightarrow {\mathbf a}\simeq {\mathbf b})]$}
%\end{enumerate}
\paragraph{III. Classification Axioms}

\begin{itemize}
\item[] {\bf Elementary Comprehension Axiom (ECA)}
\item[] $\exists X[X \simeq \{x:\psi( x)\} \wedge \forall x (x
     \vin X \leftrightarrow \psi( x))]$
\item[]
for each elementary formula  $\psi a$, which may contain additional
     parameters.
\item[]{\bf Natural Numbers}
     \begin{enumerate}
\item[(i)] ${\mathbf 0}\vin
{\mathbf N}\wedge \forall x(x \vin {\mathbf N}\rightarrow x' \vin {\mathbf N})$
       \item[(ii)] $\phi(\mathbf 0) \wedge
\forall x(\phi( x) \rightarrow \phi( x')) \rightarrow (\forall x \vin {\mathbf N})
\phi( x)$
         for each formula  $\phi$  of ${\lt}(\TO^i)$.
     \end{enumerate}
\item[]  {\bf Join (J)}
\item[]$\forall x\vin A\,
\exists Y fx\simeq Y \rightarrow \exists X[X \simeq \bj(A,f)\wedge
     \forall z(z \vin X
\leftrightarrow \exists x{\vin} A\exists y(z \simeq (x,y) \wedge y \vin
     fx))]$
\item[]{\bf Inductive Generation (IG)}
\item[]$\exists X[X \simeq \bi(A,B)
\wedge \forall x \vin A[\forall
     y[(y,x) \vin B \rightarrow y \vin X]
     \rightarrow x \vin X] \\[1ex]\phantom{XXX}\wedge\;
[\forall x
\vin A\,[\forall y \,((y,x) \vin B \rightarrow \phi (y))
\rightarrow \phi (x)]     \rightarrow \forall x \vin X\,\phi (x)]]$\\[2ex]
     where  $\phi$  is an arbitrary formula of $\TO^i$. \end{itemize}

%The first group of theories are type theories. This will be followed by constructive
\section{Type theories} The type theory   of Martin-L\"of from the 1984 book \cite{ml84} will be notated by $\MLTT^{\ext}$ where the superscript is meant to convey that this is an {\em extensional}  theory. It has all the usual type constructors $\Pi,\Sigma,+,{\mathbf 0},\mathbf{1},\mathbf{2},\mathsf{Id},\mathsf{W}$  for dependent products, dependent sums,
disjoint unions, empty type, unit type, Booleans, propositional identity types, and $\mathsf{W}$-types, respectively.
Moreover, the system comprises a sequence of universe types $\mathcal{U}_0,\mathcal{U}_1,\mathcal{U}_2,\ldots$
externally indexed by the natural numbers. The universe types are closed under the type constructors from the first list
and they form a cumulative hierarchy in that $\mathcal{U}_n$ is a type in $\mathcal{U}_{n+1}$ and if $A$ is a type in $\mathcal{U}_{n}$ then $A$ is also a type in $\mathcal{U}_{n+1}$.

In the version of \cite{ml84} the identity type was taken to be extensional whereas in the more recent versions, e.g. \cite{NPS} and the one forming the basis for homotopy type theory (see  \cite{hott}), it is considered to be intensional.
 The intensional version will simply be denoted by $\MLTT$.
  For the proof-theoretic strength, though, it turns out that the difference  is immaterial. The
  reasons will be explained below, but  perhaps a first good approximation comes from the observation
 that (exact) lower bounds can be established by interpreting certain set theories in type theory in such a way that the extensional identity type can be dispensed with in these interpretations, although for validating certain forms of the axiom of choice, e.g. the ${\mathbf{\Pi\Sigma}}\mathsf{W}\mbox{-}\AC$ axiom to be discussed below, chunks of extensionality are still required. Since we shall be discussing (partial) conservativity results of extensional over intensional type theory below,
 let's recall the differences.

 \begin{deff}{\em  A key feature of Martin-L\"of's type theory is the distinction of two notions of identity (or equality).
 {\em Judgemental identity}  appears  in judgements in the two forms $\Gamma\vdash s=t\intt A$
 and $\Gamma\vdash A=B \mbox{ type}$ between terms and between types, respectively. The general equality rules (reflexivity, symmetry, transitivity) and substitution
 rules, simultaneously at the level of terms and types, apply to these judgements as further inference rules.\footnote{See \cite[Ch.5]{NPS} or \cite[A.2.2]{hott}, where they are called structural rules.} But there is also {\em propositional identity} which gives rise to types $\bI(A,s,t)$ and allows for internal reasoning about identity.

The %formation, introduction, uniqueness and reflection
rules for the
extensional identity type are the following:\footnote{The rules are essentially the ones used in \cite{ml84}, except that
\cite{ml84} has a constant $\mathsf{r}$ as the sole canonical
element of all inhabited types $\bI(A,a,b)$. Here we use $\refl(a)$ to make the comparison with the intensional case
more transparent. In \cite{ml84}, $\bI\mbox{--Uniqueness}$ and $\bI\mbox{--Reflection}$ are called I-equality and I-elimination, respectively.}

\begin{eqnarray*}
(\bI\mbox{--Formation})&&\quad{\Gamma\vdash A\mbox{ type}\qquad \Gamma\vdash a\intt A\qquad \Gamma\vdash b\intt A
\over \Gamma\vdash\bI(A,a,b)\mbox{ type}}
\\ &&
\vphantom{ffffffffff}\\
(\bI\mbox{--Introduction})&&\quad{\Gamma\vdash a\intt A\over
\Gamma\vdash\refl(a)\intt\bI(A,a,a)}\\ \vphantom{ffffffffff} && \\
(\bI\mbox{--Uniqueness})&&\quad{\Gamma\vdash p\intt \bI(A,a,b)\over
\Gamma\vdash p=\refl(a)\intt\bI(A,a,b)}\\ \vphantom{ffffffffff} && \\
(\bI\mbox{--Reflection})&& \quad{\Gamma\vdash p\intt\bI(A,a,b)\over \Gamma\vdash  a=b\intt
A}.\end{eqnarray*}
Reflection has the effect of rendering judgemental identity undecidable, i.e., the (type checking) questions whether $\Gamma\vdash a=b\intt A$  or $\Gamma\vdash a\intt A$ hold  become undecidable.
On the other hand, the set-theoretic models and many recursion-theoretic  models of type theory (see \cite{relating,beeson,rohio}) validate extensionality, lending it an intuitive appeal.

 For the intensional identity type, the foregoing rules of formation and introduction are retained, however, uniqueness and reflection are jettisoned,
 getting replaced by elimination and equality rules which are motivated by Leibniz's principle of indiscernibility, namely that identical elements are those that satisfy the same properties. Though instead of capturing identity by quantifying (impredicatively) over all properties (as in Principia),
 the entire family of identity types $(\bI(A,x,y))_{x,y\intt A}$ %over $A\times A$
 is viewed as being inductively generated with sole constructor $\refl$ (see \cite{NPS,hott}).
 The elimination and equality rules are the following:
\begin{eqnarray*}
(\Idd\mbox{--Elimination})&&\quad{\dekordfuenf{\Gamma\vdash a\intt A}{\Gamma\vdash b\intt
A}{\Gamma\vdash c\intt \Idd(A,a,b)}{\Gamma,\,x\intt A,\,y\intt A,\, z\intt\Idd(A,x,y)\vdash C(x,y,z)\mbox{ type}} {\Gamma,\,x\intt
A\vdash d(x)\intt C(x,x,\refl(x))} \over \Gamma\vdash \idpeelJ(c,d)\intt C(a,b,c)}\\ \vphantom{ffffffffff} &&\\
(\Idd\mbox{--Equality})&&\quad{\dekorddrei{\Gamma\vdash a\intt
A}{\Gamma,\,x\intt A,\, y\intt A,\,
z\intt\Idd(A,x,y)\vdash C(x,y,z)\mbox{ type}} {\Gamma,\,x\intt A\vdash d(x)\intt C(x,x,\refl(x))}
\over \Gamma\vdash  \idpeelJ(\refl(a),d)=d(a)\intt C(a,a,\refl(a))\,.}
\end{eqnarray*}
}\end{deff}
An immediate consequence of these rules is the indiscernibility of identical elements expressed as follows. For every family $(C(x))_{x\intt A}$ of types there is a function
$$f\intt \Pi_{x,y\intt A}\Pi_{p\intt {\mathsf{Id}}(A,x,y)}[C(x)\to C(y)]$$
such that with $1_{C(x)}$ being the function $u\mapsto u$ on $C(x)$ we have $f(x,x,\refl(x))= 1_{C(x)}$.

Foregoing extensional identity and using
the induction principle encapsulated in $\Idd\mbox{-elimination}$ and $\Idd\mbox{-equality}$ in its stead, is crucial to the more subtle homotopy interpretations of type theory.

\section{Constructive set theories}
Constructive Set Theory was introduced by
 Myhill in a seminal paper \cite{myhill}, where a
specific axiom system $\CST$ was introduced.
 Through developing constructive
 set theory he wanted to isolate the principles underlying Bishop's conception
of what sets and functions are, and he  wanted  ``these
principles to be such as to
 make the process of formalization completely trivial,
 as it is in the classical case"
(\cite{myhill}, p. 347).
 Myhill's $\CST$ was subsequently modified by Aczel
 and  the resulting
  theory was called {\em Constructive Zermelo-Fraenkel set theory},
  $\CZF$. A hallmark of this theory is that it possesses a type-theoretic
interpretation (cf. \cite{aczel82,mar}).
Specifically, $\CZF$ has a scheme called Subset Collection Axiom
(which is a generalization of Myhill's Exponentiation Axiom) whose
formalization was directly inspired by the type-theoretic
interpretation.

The language of $\CZF$ is the same first order language
as that of classical Zermelo-Fraenkel Set Theory, $\ZF$ whose only
non-logical symbol is $\in$. The logic of $\CZF$ is intuitionistic
first order logic with equality. Among its non-logical axioms are
{\it Extensionality}, {\it Pairing} and {\it Union} in their usual
forms.
  $\CZF$ has
additionally axiom schemata which we will now proceed to
summarize. Below $\emptyset$ stands for the empty set and $v+1$ denotes $v\cup\{v\}$.
 A set-theoretic formula
is said to be {\em restricted} or {\em bounded} or $\Delta_0$ if it is constructed from prime formulae using
$\neg,\wedge,\vee,\rightarrow$ and only restricted quantifiers $\forall x\ini y,\; \exists x\ini
y$.
\\[1ex] {\em Infinity:}\footnote{This axiom asserts the existence of a unique set usually called $\omega$. Note that the second conjunct
in $[\ldots]$ entails the usual induction principle for $\omega$ with regard to set properties (or equivalently $\Delta_0$ formulae).}
$$\exists x\,[\forall u\bigl(u\ini x\leftrightarrow
\bigl(\emptyset=u\oder \exists v\ini x\;u=v+1 \bigr)\bigr)\;\wedge\;\forall z\,(\emptyset\in z\,\wedge\, \forall y\in z\;y+1\in z\to x\subseteq z)].$$
\\[1ex] {\em Set Induction:} For all formulae $\phi$,
$$\forall x[\forall y \in x \phi (y) \rightarrow \phi (x)]
\rightarrow \forall x \phi (x).$$
{\em Restricted or Bounded Separation:} For all {\it restricted} formulae  $\phi$,
$$\forall a \exists b \forall x
[x \in b \leftrightarrow x \in a \wedge \phi (x)].$$
{\em Strong Collection:}  For all formulae $\phi$,
$$\forall a
\bigl[\forall x \in a \exists y \phi (x,y)\;\rightarrow\; \exists
b \,[\forall x \in a\, \exists y \in b\, \phi (x,y) \wedge \forall
y \in b\, \exists x \in a \,\phi (x,y)]\bigr].$$
{\em Subset Collection:}  For all formulae $\psi$,
\begin{eqnarray*}\lefteqn{\forall a \forall b \exists c
\forall u \,\bigl[\forall x \in a \,\exists y \in b\;\psi (x,y,u)
\,\rightarrow \,}\\ && \exists d \in c \,[\forall x \in a\,
\exists y \in d\, \psi (x,y,u) \wedge \forall y \in d\, \exists x
\in a \,\psi (x,y,u)]\bigr].\end{eqnarray*}
  The Subset Collection
schema easily qualifies as the most intricate axiom of $\CZF$.

We shall also consider an additional axiom that holds true in the type-theoretic interpretation of Aczel
if the type theory is equipped with $W$-types. To introduce it, we need the notion of a regular set.
The formula in the language
of \czf\ defining the property of a set $A$ that it is {\it regular} states
that $A$ %the set in question
is transitive, and for every $a\in A$ and set $R\subseteq a\times A$ if
$\forall x\in a\,\exists y\,( \langle x,y\rangle\in R)$,
then there is a set $b\in A$
such that
$$\forall x\in a\,\exists y\in b\,( \langle x,y\rangle\in R)\;\wedge\;
\forall y\in b\,\exists x\in a\, (\langle x,y\rangle\in  R).$$
In particular,
 if  $R:\,a\to A$ is a function, then the image of $R$ is an
an element of $A$.
%we have constructed a function from an element of
%that set into the set in
%questions, then there is an element of it such that the image of the
%function is contained in it.
Let  $\mathsf{Reg}(A)$  denote this assertion.  With this auxiliary
definition we
can state the
\\[1ex]
{\em Regular Extension Axiom $\REA$}
$$\forall x \exists y [x \subseteq y \wedge \mathsf{Reg}(y)]\,.$$

%There are several yet stronger forms of the axiom of choice that are validated in the type-theoretic interpretation
%(see \cite{aczel82}). Stating them requires some preliminary definitions.

\subsection{The axiom of choice in constructive set theories}
Among the axioms of set theory, the axiom of choice is
distinguished by the fact that  it is the only one that one finds
mentioned in  workaday mathematics. In the mathematical world of
the beginning of the 20th century, discussions about the status of
the axiom of choice were important. In 1904 Zermelo proved that
every set can be well-ordered by employing the axiom of choice.
While Zermelo argued that it was self-evident, it was also
criticized as an excessively non-constructive principle by some of
the most distinguished analysts of the day, notably Borel, Baire,
and Lebesgue.
  At first blush this reaction
against the axiom of choice utilized in Cantor's new theory of
sets is surprising as the French analysts had used and continued
to use choice principles routinely in their work. However, in the
context of 19th century classical analysis only the Axiom of
Dependent Choices, $\DC$, is invoked and considered to be natural,
while the full axiom of choice is unnecessary and even has some
counterintuitive consequences.

Unsurprisingly, the axiom of choice does not have a unambiguous
status  in constructive mathematics either. On the one hand it is
said to be an immediate consequence of the constructive
interpretation of the quantifiers. Any proof of $\forall x\inn
A\,\exists y\inn B\,\phi(x,y)$ must yield a function $f:A\ra B$
such that $\forall x\inn A\,\phi(x,f(x))$. This is certainly the
case in Martin-L\"of's intuitionistic theory of types. On the
other hand,
 it has been observed that the full
axiom of choice cannot be added to systems of extensional
constructive set theory without yielding constructively
unacceptable cases of excluded middle (see \cite{Diaconescu}). In
extensional intuitionistic set theories, a proof of a statement
$\forall x\inn A\,\exists y\inn B\,\phi(x,y)$, in general,
provides only a function $F$, which when fed a proof $p$
witnessing  $x\inn A$, yields $F(p)\inn B$ and $\phi(x,F(p))$.
Therefore, in the main, such an $F$ cannot be rendered a function
of $x$ alone. Choice will then hold over sets which have a
canonical proof function, where a constructive function $h$ is a
canonical proof function for $A$ if for each $x\inn A$, $h(x)$ is
a constructive proof that $x\inn A$. Such sets  having natural
canonical proof functions ``built-in" have been called {\em bases}
(cf. \cite{TD88}, p. 841).

\paragraph{Some constructive choice principles} In  many a text on
constructive mathematics,  axioms of countable choice and
dependent choices are accepted as constructive principles. This
is, for instance, the case in Bishop's constructive mathematics
(cf. \cite{bishop}) as well as Brouwer's intuitionistic analysis
(cf. \cite{TD88}, Ch. 4, Sect. 2). Myhill also incorporated  these
axioms in his constructive set theory \cite{myhill}.

The weakest constructive choice principle we shall consider is the
{\em Axiom of Countable Choice}, $\ACo$, i.e. whenever $F$ is a
function with   domain $\omega$ such that $\forall i\inn
\omega\,\exists y\inn F(i)$, then there exists a function $f$ with
domain $\omega$ such that $\forall i\inn \omega\,f(i) \inn F(i)$.

A  mathematically very useful axiom to have in set theory is the
{\em Dependent Choices Axiom}, $\DC$, i.e., for all formulae
$\psi$, whenever
$$(\forall x\inn a)\,(\exists y\inn a)\,\psi(x,y)$$ and $b_0\inn a$, then
there exists a function $f:\omega\rightarrow a$ such that
$f(0)=b_0$ and
$$(\forall n\inn\omega)\,\psi(f(n),f(n+1)).$$
Even more useful is the
 {\em Relativized Dependent Choices Axiom}, $\RDC$.
It asserts that for arbitrary formulae $\phi$ and $\psi$, whenever
$$\forall x\bigl[\phi(x)\,\rightarrow\,\exists y
\bigl(\phi(y)\,\wedge\,\psi(x,y)\bigr)\bigr]$$ and $\phi(b_0)$,
then there exists a function $f$ with domain $\omega$ such that
$f(0)=b_0$ and
$$(\forall n\inn\omega)\bigl[\phi(f(n))\,\wedge\,\psi(f(n),f(n+1))\bigr].$$

In addition to the ``traditional" axioms  of choice stated above, the interpretation of
set theory in type theory  validates several  new choice principles which are are not well known. To state them we need to
introduce various operations on classes.

\begin{rema}{\em Let $\CZFe$ denote the modification of $\CZF$
with Eponentiation in place of Subset
Collection.

 In almost all the results of this paper, $\CZF$ could be replaced by $\CZFe$, that is to say, for
the purposes of this paper it is enough to assume Exponentiation
rather than Subset Collection. However, in what follows we shall
not point this out again.}\end{rema}

\begin{definition} \label{dePiSi} $(\CZF)$
{\em If $A$ is a set and $B_x$ are classes for all $x\in
A$, we define a class $\prod_{x\in A}B_x$ by:  \begin{eqnarray}\label{PiSi1}
\prod_{x\in A}B_x &:=& \{f\mid f:A \to {\bigcup_{x\in A}B_x} \,\wedge\,\A
x\lin A (f(x)\in B_x)\}. \end{eqnarray} If $A$ is a class and $B_x$ are
classes for all $x\in A$, we define a class $\sum_{x\in A}B_x$ by:
 \begin{eqnarray}\label{PiSi2} \sum_{x\in A}B_x &:=& \set {\pair{x}{y}}
{x\in A \land y\in B_x}. \end{eqnarray} If $A$ is a class and $a,b$ are
sets, we define a class $\I(A,a,b)$ by:  \begin{eqnarray} \label{PiSi3}
\I(A,a,b)&:=& \set {z\in 1} {a=b\,\wedge\,a,b\in A}. \end{eqnarray} If $A$
is a class and for each $a\in A$, $B_a$ is a set, then
$${\mathbf W}_{a\in A}B_a$$ is the smallest class $Y$ such that whenever $a\in A$ and
$f:B_a\rightarrow Y$, then $\langle a,f\rangle\in Y$. }\end{definition}

\begin{lemma}\label{lePiSi} $(\CZF)$
 If $A$,$B$,$a$,$b$ are sets and $B_x$ is a set for
all $x\in A$, then $\prod_{x\in A}B_x$, $\sum_{x\in A}B_x$ and
$\I(A,a,b)$ are sets. \end{lemma} \Proof\ \cite[Lemma 2.5]{rt}.\qed

In the following we shall introduce several inductively defined classes, and, moreover,
we have to ensure that such classes
can be formalized in $\CZF$.

We define an {\em inductive definition} to be a class of ordered pairs.
If $\Phi$ is an inductive definition and $\langle x,a\rangle\in\Phi$ then
we write
    \begin{eqnarray*} &&\step{x}{a}{_\Phi}\end{eqnarray*}
   %$$ \frac{x}{a}{_\Phi}$$
and call $\step{x}{a}{}$ an {\em (inference) step} of $\Phi$, with set
$x$ of {\em
premisses} and {\em conclusion} $a$.  For any class $Y$, let
\begin{eqnarray*} \Gamma_{\Phi}(Y)&=& \bigl\{a\,\mid\;\exists x\,\bigl(x\subseteq Y\;\;\wedge\;\;
\step{x}{a}{_\Phi}\,\bigr)\bigr\}.\end{eqnarray*}
 The class $Y$ is {\em $\Phi$-closed} if $\Gamma_{\Phi}(Y)\subseteq Y$.
Note that $\Gamma$ is monotone; i.e. for classes $Y_1,Y_2$, whenever $Y_1\subseteq Y_2$, then
$\Gamma(Y_1)\subseteq \Gamma(Y_2)$.

We define the class {\em inductively defined by $\Phi$} to be the
smallest $\Phi$-closed class.  The main result about inductively defined classes states
that this class, denoted $\I(\Phi)$, always exists.
\ble\label{ind} $(\CZF)$ {\em (Class Inductive Definition Theorem)}
For any inductive definition $\Phi$ there is a smallest $\Phi$-closed
class $\I(\Phi)$.
\ele
\Proof\
\cite{aczel82}, section 4.2 or \cite{maralt}, Theorem 5.1. \QED

\ble\label{lePiSiW} $(\CZF+\REA)$ If $A$ is a  set and $B_x$ is a
set for all $x\in A$, then ${\mathbf W}_{a\in A}B_a$ is a set.
\ele \Proof\ This follows from \cite{aczel86}, Corollary 5.3. \QED

\ble\ $(\CZF)$\\
\label{leY} There exists a smallest $\mathbf{\Pi\Si}$-closed
class, i.e., a smallest class \BY\ such that the following
hold:\\
$(i)\;$ $n\in\BY$ for all $n\in\omega$;\\
$(ii)\;$ $\om\in\BY$;\\
$(iii)\;$ $\prod_{x\in A}B_x\in\BY$ and $\sum_{x\in A}B_x\in\BY$
whenever $A\in\BY$ and $B_x\in\BY$ for all $x\in A$.

 Likewise,
there exists a smallest $\mathbf{\Pi\Si}\I$-closed class, i.e. a
smallest class \BYst, which, in addition to the closure conditions
$(i)$--$(iii)$ above,
satisfies:\\[1ex]
%Hier
$(iv)\;$ $\I(A,a,b)\in\BYst$ whenever $A\in\BYst$ and $a,b\in A$.
\ele \Proof\  \cite[Lemma 2.8]{rt}. \qed

\begin{deff}{\em The $\mathbf{\Pi\Sigma}$-generated sets are the sets in the smallest
$\mathbf{\Pi\Sigma}$-closed class. Similarly one
defines the $\mathbf{\Pi\Sigma} \I$, $\mathbf{\Pi\Sigma} \WW$ and
$\mathbf{\Pi\Sigma} \WW\I$-generated sets.

A set $P$ is a {\em base} if
for any $P$-indexed family $(X_a)_{a\in P}$ of inhabited sets $X_a$, there exists a function
$f$ with domain $P$ such that, for all $a\in P$, $f(a)\in X_a$.

\PSAC\ is the statement that every $\mathbf{\Pi\Sigma}$-generated set is a base. Similarly
one defines the axioms \PSIAC, \PSWIAC, and \PSWAC.

The {\em presentation axiom}, $\PAx$, states that every set is the surjective image of a base. }\end{deff}

 \begin{lem} \label{coYY*} \
\begin{itemize}
\item[(i)] $(\CZF)$
\PSAC\ and \PSIAC\ are equivalent.
\item[(ii)] $(\CZF+\REA)$
\PSWAC\ and \PSWIAC\ are equivalent.
\end{itemize}
\end{lem}
 \Proof\ \cite[2.12]{rt}. \qed

\subsection{Large sets in constructive set theory}
Large cardinals play a central role in modern set theory. This
section deals with large cardinal properties in the context of
intuitionistic set theories. Since in intuitionistic set theory
$\in$ is not a linear ordering on ordinals the notion of a
cardinal does not play a central role. Consequently, one talks
about {\em ``\,large set properties''} instead of {\em ``\,large
cardinal properties''}. %Friedman and \v{S}\v{c}edrov \cite{fsc}
%studied large set properties in the context of $\IZF$.
When stating these properties one has to proceed rather carefully.
Classical equivalences of cardinal notion might no longer prevail
in the intuitionistic setting, and one therefore wants to choose a
rendering which intuitionistically retains the most strength. On
the other hand certain notions have to be avoided so as not to
imply excluded third. To give an example, cardinal notions like
measurability, supercompactness and hugeness have to be expressed
in terms of elementary embeddings rather than ultrafilters.

We shall, however, not concern ourselves with very large cardinals
here and rather restrict attention to the very first notions of
largeness introduced by Hausdorff and Mahlo, that is, inaccessible
and Mahlo sets and the pertaining hierarchies of inaccessible and
Mahlo sets.

We have already seen one notion of largeness, namely that of a regular set.
In $\ZFC$, a regular set which itself is a model of the axioms of $\CZF$ is of the form $V_{\kappa}$ with $\kappa$ a strongly inaccessible cardinal.\footnote{Note that $\CZF$ with classical logic is the same theory as $\ZF$.} In the context of $\CZF$ this notion is much weaker.

\begin{deff}\label{4.4}{\em
If $A$ is a transitive set and $\phi$ is a formula with parameters
in $A$ we denote by $\phi^A$ the formula which arises from $\phi$
by replacing all unbounded quantifiers $\forall u$ and $\exists v$
in $\phi$ by $\forall u\in A$ and $\exists v\in A$, respectively.

We can view any transitive set $A$ as a structure equipped with
the binary relation ${\in}_A=\{\paar{x,y}\mid x\in y\in A\}$. A
set-theoretic sentence whose parameters lie in $A$, then has a
canonical interpretation in $(A,{\in}_A)$ by interpreting $\in$ as
${\in}_A$, and $(A,{\in}_A)\models \phi$ is logically equivalent
to $\phi^A$. We shall usually  write $A\models \phi$ in place of
$\phi^A$.

A set $I$ is said to be {\em weakly inaccessible} if $I$ is a regular set
such that $I\models \CZF^-$, where $\CZF^-$ denotes the theory $\CZF$ bereft of the set induction scheme.\footnote{
Note that if the background set theory
validates set induction for $\Delta_0$ formulae then a transitive set will be automatically a model of the full set induction scheme, and thus
a regular set $I$ will satisfy $I\models \CZF$.}

The strong regular extension axiom, $\sREA$, states that every set is an element of  a weakly
inaccessible set.
}
\end{deff}

There is a more `algebraic' way of expressing weak inaccessibility.
Stating it requires some definitions.

\begin{deff}{\em For sets $A,B$ we denote by $\voll{A}{B}$ the collection of all full relations from $A$ to $B$, i.e.,
of those relations $R\subseteq A\times B$ such that $\forall x\in A\,\exists y\in B\; \langle x,y\rangle \in R$.
A set $C$ is said to be {\em full in $\voll{A}{B}$} if for all $R\in \voll{A}{B}$ there exists $R'\in\voll{A}{B}$
such that $R'\subseteq R$ and $R'\in C$.

For a set $A$, define $\bigwedge A$ to be the set $\{x\in 1\mid \forall u\in A\,x\in u\}$, where $1=\{\emptyset\}$.
}\end{deff}

\begin{prp} $(\CZF^-)$
A set $I$ is  weakly inaccessible if and only if  $I$ is a regular set
such that the following are satisfied:
 \begin{enumerate}
\item $\omega\in I$,
 \item $\forall a\in I\;\bigcup a\in I$,
 \item $\forall a\in I\,[a\mbox{ inhabited}\;\rightarrow \;\bigcap
 a\in I]$,
  \item $\forall A,B\in I\,\exists C\in I\;\;\mbox{
 $C$ is full in $\voll{A}{B}$}$.
 \end{enumerate}
 \end{prp}
 \prf \cite[10.26]{mar}.

 We will consider two stronger notions.

 \begin{deff}\label{4.7}{\em A set $I$ is called {\em inaccessible} if $I$ is weakly inaccessible
 and for all $x\in I$ there exists a regular set $y\in I$ such that $x\in y$.

 A set $M$ is said to be {\em Mahlo} if $M$ is
inaccessible
and for every $R\in\voll{M}{M}$ there exists an inaccessible $I\in M$
such that $$\forall x\in I\,\exists y\in I\;\langle x,y\rangle\in R.$$
}\end{deff}

 \subsection{Fragments of second order arithmetic}
 The proof-theoretic strength of theories is commonly calibrated
using standard theories and their canonical fragments. In
classical set theory this linear line of consistency strengths is
couched in terms of large cardinal axioms while for weaker
theories the line of reference systems traditionally consist of
subsystems of second order arithmetic. The observation that large chunks of mathematics
can already be formalized in fragments of second order arithmetic
goes back to Hilbert and
Bernays \cite{hb}, and has led to a systematic research program known as {\em Reverse Mathematics}.
Below we give an  account of the  syntax of ${\mathcal L}_2$ and frequently considered axiomatic principles.
\begin{dfn}{\em
The language ${\mathcal L}_2$ of second-order arithmetic contains
 number variables $x,y,z,u,\ldots$, set variables
$X,Y,Z,U,V,A,B,C,\ldots$ (ranging over subsets of $\mathbb N$),
the constant $0$, function symbols $Suc,+,\cdot$, and relation
symbols $=,<,\in$. $Suc$ stands for the successor function. We write $x+1$ for $Suc(x)$. {\em
Terms} are built up as usual. For $n\!\in\!\mbox{$\mathbb N$}$,
let $\bar{n}$ be the canonical term denoting $n$. Formulae are
built from the prime formulae $s=t$, $s<t$, and $s\in X$ using
$\wedge,\vee,\neg,\forall x,\exists x, \forall X$ and $\exists X$
where $s,t$ are terms. Note that equality in ${\mathcal L}_2$ is
only a relation on numbers. However, equality of sets will be
considered a defined notion, namely $X=Y$ if and only if $\forall
x[x\!\in\! X\,\leftrightarrow\, x\!\in\! Y]$. As per usual, number
quantifiers are called bounded if they occur in the context
$\forall x(x<s\rightarrow\ldots)$ or $\exists x (x<s\wedge\ldots)$
for a term $s$ which does not contain $x$. The
${\mathbf{\Sigma}}^0_0$-formulae are those formulae in which all
quantifiers are bounded number quantifiers. For $k>0$,
${\mathbf{\Sigma}}^0_k$-formulae are formulae of the form $\exists
x_1\forall x_2\ldots Qx_k\phi$, where $\phi$ is
${\mathbf{\Sigma}}^0_0$; ${\mathbf{\Pi}}^0_k$-formulae are those
of the form $\forall x_1\exists x_2\ldots Qx_k\phi$. The union of
all ${\mathbf{\Pi}}^0_k$- and ${\mathbf{\Sigma}}^0_k$-formulae for
all $k\in{\mathbb N}$ is the class of {\em arithmetical} or {\em
${\mathbf{\Pi}}^0_\infty$-formulae}. The
${\mathbf{\Sigma}}^1_k$-formulae (${\mathbf{\Pi}}^1_k$-formulae)
are the formulae $\exists X_1\forall X_2\ldots QX_k\phi$
(resp.~$\forall X_1\exists X_2\ldots Qx_k\phi$) for arithmetical
$\phi$.

The basic axioms in all theories of second-order arithmetic are
the defining axioms of $0,1,+,\cdot,<$ and the {\em induction
axiom}
\[
\forall X(0\in X\wedge\forall x(x\in X\rightarrow x+1\in
X)\rightarrow
          \forall x(x\in X)),
\]
respectively the {\em scheme of induction}
\[ {\mathbf{IND}}\qquad
\phi(0)\wedge\forall
x(\phi(x)\rightarrow\phi(x+1))\rightarrow\forall x\phi(x),
\]
where $\phi$ is an arbitrary ${\mathcal L}_2$-formula.
We consider
the axiom scheme of {\em ${\mathcal C}$-comprehension}  for formula
classes ${\mathcal C}$ which is given by
\[{\mathcal C}\mbox{-}{\mathbf{CA}}\qquad\exists X\forall u(u\in X\leftrightarrow\phi(u))
\]
 for all formulae $\phi\in{\mathcal C}$ (of course, $X$ must not be free in $\phi$).

 For each axiom scheme $\mathbf{Ax}$ we denote by $(\mathbf{Ax})$
 the theory consisting
of the basic arithmetical axioms, the scheme
${\mathbf{\Pi}}^0_\infty\mbox{-}\mathbf{CA}$, the scheme of
induction and the scheme $\mathbf{Ax}$. If we replace the scheme
of induction by the induction axiom, we denote the resulting
theory by $(\mathbf{Ax})_0$. An example for these notations is the
theory $({\mathbf{\Pi}}^1_1\mbox{-}{\mathbf{CA}})$ which contains
the induction scheme, whereas
$({\mathbf{\Pi}}^1_1\mbox{-}{\mathbf{CA}})_0$ only contains the
induction axiom in addition to the comprehension scheme for
$\mathbf \Pi^1_1$-formulae.

In the basic system  one can introduce defined
symbols for all primitive recursive functions. Especially, let
$\langle{,}\rangle:{\mathbb N}\times{\mathbb
N}\longrightarrow{\mathbb N}$ be a primitive recursive and
bijective pairing function. The $x^{th}$ section of $U$ is defined
by $U_x\,:=\,\{y:\,\langle{x,y}\rangle\in U\}$. Observe that a set
$U$ is uniquely determined by its sections on account of
$\langle{,}\rangle$'s bijectivity. Any set $R$ gives rise to a
binary relation $\prec_{R}$ defined by
$y\prec_{R}x\,:=\,\langle{y,x}\rangle\in R$.
Using  this coding we can formulate the
 {\em ${\mathcal C}$-axiom of choice} scheme for formula
classes ${\mathcal C}$ which is given by
\[{\mathcal C}\mbox{-}{\mathbf{AC}}\qquad\forall x\exists Y\,\psi(x,Y) \to \exists Z\,\forall u\,\psi(x,Z_x),
\]
for all formulae $\psi\in{\mathcal C}$ ($Z$ must not be free in $\psi$).

Another important principle is {\em Bar induction}:
$${\mathbf{BI}}\qquad\forall X\bigl[
 {\mathbf{WF}}(\prec_X)\;\wedge\;\forall u\bigl(\forall v\prec_X u\phi(v)\rightarrow
\phi(u)\bigr) \,\rightarrow\, \forall u\phi(u)\bigr]$$ for all
formulae $\phi$, where ${\mathbf{WF}}(\prec_X)$ expresses that
$\prec_X$ is well-founded, i.e., ${\mathbf{WF}}(\prec_X)$ stands
for the formula $$\forall Y\,\bigl[\forall u\bigl[(\forall v\prec_X
u\;v\in Y)\,\rightarrow\, u\in Y\bigr] \,\rightarrow\, \forall
u\;u\in Y\bigr].$$
 }\end{dfn}
 Universes in type theory (with $W$-types) bear a strong relation to $\beta$-models which are models of the language of
 ${\mathcal L}_2$ or set theory for which the notion well-foundedness is absolute.
 \begin{deff}\label{4.8}{\em Any set $A$ of natural numbers  gives rise to a set $\mathfrak{X}_A:=\{A_i\mid i\in\mathbb{N}\}$ of sets of natural numbers. $A$ is  said to be a {\em $\beta$-model} if the ${\mathcal L}_2$-structure $$\mathfrak{A}:=(\mathbb{N},\mathfrak{X}_A, 0,1,+,\cdot,\in)$$
 is a $\beta$-model, i.e., $\mathfrak{A}\models {\mathbf{\Pi}}^0_\infty\mbox{-}\mathbf{CA}$, and
  whenever $Y\in \mathfrak{X}_A$ and $\mathfrak{A}\models {\mathbf{WF}}(\prec_Y)$ then
  $\prec_Y$ is well-founded.

  Obviously, the notion, the notion of $\beta$-model can be expressed in ${\mathcal L}_2$.}\end{deff}

  \paragraph{An intuitionistic ${\mathcal L}_2$-theory.} There is an interesting version of  second order arithmetic,
 which will be used in theory reductions,  that classically has the same
  strength as full second order arithmetic, $(\Pi^1_{\infty}\mbox{-}\mathbf{CA})$, but when based on intuitionistic logic is of the same strength
  as $\TO^i$.

\begin{dfn}\label{4.2a}{\em %$(The theory \Iario)
\Iari is a theory in the language of second order arithmetic.
  The
logical rules of \Iari are those of intuitionistic second order
arithmetic.
In addition to the usual axioms for intuitionistic second
order logic, axioms are (the universal closures of):
\begin{enumerate}
\item {\bf Induction}:
$$\phi(0)\land\forall n[\phi(n)\then\phi(n+1)]\then\forall n\phi(n)$$
for all formulae $\phi$.
\item {\bf Arithmetic Comprehension Schema}:
$$\exists X\forall n[n\in X\lra\psi(x)]$$
for $\psi$ arithmetical (parameters allowed).
\item {\bf Replacement}: %Axiom of Countable Collection:
$$\forall X[\forall n\in X\exists\,!
Y\phi(n,Y)\then\exists Z\forall n\in X\,\phi(n,
 Z_n)]$$
for all formulas $\phi$. Here $\phi(n,Z_n)$ arises from $\phi(n,Z)$ by
replacing each occurrence $t\in Z$ in the formula by $\pair nt\in Z$.
\item {\bf Inductive Generation}:
$$\forall U\forall X\exists Y\bigl[
\WPU(X,Y)\land(\forall n[\forall k(k\ltx n\then\phi(k))\then\phi(n)]\then
\forall m\in Y\,\phi(m))\bigr],$$
for all formulas $\phi$, where $k\ltx n$ abbreviates $\pair kn\in X$ and
$\WPU(X,Y)$ stands for
$$\progU(X,Y)\land\forall Z[\progU(X,Z)\then Y\subseteq Z]$$
with $\progU(X,Y)$ being $\forall n \in U[\forall k(k\ltx n\then k\in Y)
\then n\in Y]$.
\end{enumerate}}
\end{dfn}

\begin{rmk}\label{4.3a}{\em (\Iario)
Note that {$\WPU(X,Y)$} and {$ \WPU(X,Y')$} imply $Y=Y'$,
i.e. $\forall n(
n\in Y\lra n\in Y')$. Therefore, if {$\WPU(X,Y)$}, then
$$\forall n\in U[\forall k\ltx n\,\phi(k)\then\phi(n)]\then\forall m\in Y\,
\phi(m)$$
holds for all formulae $\phi$.}
\end{rmk}

The latter principle will be referred to as {\em `` induction over the
well--founded part of $\ltx$'' }. In the rest of this section we shall write
$\wof(U,X)$ for the (extensionally) uniquely determined $Y$ which satisfies
$\WPU(X,Y)$.

The main tool for performing  the well-ordering proof of
\cite{j83} in \Iari is the following principle of transfinite recursion.
%We write $\wf(X)$ iff $\wf(X,\nn)$ holds, where
%$\nn:=\{n:n=n\}$. We want to show that tn\ can be modeled in \Iari
%similarly as in \czfp.
\begin{prp}\label{4.4a}
{\em (\Iario)} If {\em$\WPU(X,Y)$} and
$\forall n\in Y\forall W\exists!V\,\psi(n,W,V)$, then
there exists $Z$ such that $$\forall n\fin Y\,\psi(n,\bigcup\{(Z)_k:
k\ltx n\},(Z)_n).$$
\end{prp}
\prf See \cite[6.4]{rohio}.\qed

\section{On relating theories I}
The first result relates intuitionistic explicit mathematics to constructive set theory and a fragment of $\MLTT$. Let $\mathbf{MLT}_{1W}V$ be the fragment of $\MLTT$ with only one universe $\mathcal{U}_0$ where the $\mathsf{W}$-constructor can solely be applied to families of types in $\mathcal{U}_0$ but one can also form  the type $\mathsf{V}:=\mathsf{W}_{(A:\mathcal{U}_0)}A$ (something that could be called the type of Brouwer ordinals of $\mathcal{U}_0$). We shall also consider the type theory $\mathbf{MLT}_{1W}$ which is the fragment of $\mathbf{MLT}_{1W}V$ without the type $\mathsf{V}$.

\paragraph{A principle of omniscience.}
Certain basic principles of classical mathematics are taboo for
the constructive mathematician. Bishop called them {\em principles
of omniscience}.
The limited principle of omniscience, $\LPO$, is an instance of the law of excluded middle which usually
serves as a line of demarcation, separating ``constructive" from ``non-constructive" theories.
In the case of $\CZF$, adding the law of excluded middle even just for atomic statements of the form $a\in b$ results in an enormous increase in proof strength, pushing it up beyond that of Zermelo set theory. However, $\LPO$ can be added to $\CZF$ without affecting its proof-theoretic strength.
$\LPO$ has the pleasant side effect that one can carry out elementary analysis pretty much in the same way as in any standard text book.

 \begin{definition}\label{omnis} {\em Let $2^{\bN}$ be Cantor space, i.e the set of all functions from the naturals into $\{0,1\}$.
 { Limited Principle of Omniscience} ($\LPO$):
$$\forall f\in 2^{\bN}\,[\exists n\,f(n)=1\;\;\vee\;\;\forall n\,f(n)=0].$$
 }\end{definition}

\begin{thm}\label{Theorem1} The following theories have the same proof-theoretic strength and therefore prove (as a minimum)
the same $\Pi^0_2$ statements of arithmetic:
\begin{itemize} \item[(i)] Intuitionistic explicit mathematics, $\TO^i$.
\item[(ii)] Constructive Zermelo-Fraenkel set theory with the {\em Regular Extension Axiom}, $\CZF+\REA$.
\item[(iii)] Constructive Zermelo-Fraenkel set theory augmented by $\RDC$ and the {\em strong Regular Extension Axiom}, $\CZF+\sREA+\RDC$.
\item[(iv)] $\CZF+\REA+{\mathbf{\Pi\Sigma}}{\mathsf{W}}\mbox{-}\AC +\RDC + \PAx$.
\item[(v)] The extensional type theory $\mathbf{MLT}_{1W}^{\ext}\mathsf{V}$.
 \item[(vi)]   $\mathbf{MLT}_{1W}\mathsf{V}$.
\item[(vii)] The extensional type theory $\mathbf{MLT}_{1W}^{\ext}$.
\item[(viii)] $\mathbf{MLT}_{1W}$.
\item[(ix)] The classical subsystem of second order arithmetic $({\mathbf{\Sigma}}^1_2\mbox{-}\mathbf{AC})+\mathbf{BI}$
(same as $({\mathbf{\Delta}}^1_2\mbox{-}\mathbf{CA})+\mathbf{BI}$).
\item[(x)] The intuitionistic system $\mathbf{IARI}$ of second order arithmetic.
\item[(xi)] Classical Kripke-Platek set theory,$\mathbf{KP}$ (cf. \cite{barwise}, plus the axiom asserting that every set is contained in an admissible set. (This theory is often denoted by $\KPi$.)
    \item[(xii)] Intuitionistic Kripke-Platek set theory, $\mathbf{IKP}$, plus the axiom asserting that every set is contained in an admissible set. (This theory will be notated by $\mathbf{IKPi}$.)
\item[(xiii)] $\CZF+\REA+\RDC + \LPO$.

\end{itemize}
\end{thm}
\prf The equivalence of (i),(ii),(iii),(iv),(v),(vi),(vii),(viii),(ix),(x), and (xi) follows from \cite{rohio}, Theorem 3.9, Proposition 5.3, Theorem 5.13 and Theorem 6.13 plus the extra observation that the interpretation of $\mathbf{IRA}$ in $\mathbf{MLT}_{1W}^{\ext}$ defined in \cite[Definition 6.5]{rohio} and proved to be an interpretation
in \cite[Theorem 6.9]{rohio} actually only requires the intensional identity type. It was already observed by Palmgren \cite{p92} that the interpretations of theories of
iterated, strictly positive inductive definitions in type theory works with the intensional identity, and the same argument applies here.

The equivalence of (ii) and (iii) follows from \cite[Theorem 4.7]{anti}, where the principle $\sREA$ is denoted by $\mathsf{INAC}$.

The proof-theoretic equivalence of (xi) and (xii) follows since the intuitionistic version is a subtheory of the classical one and
the well-ordering proof for initial segments of the ordinal of $\mathbf{KPi}$ can already be carried out in the intuitionistic theory.

 For (xiii) we rely on \cite{lpo}. That the theory  $\CZF+\REA+\RDC + \LPO$  has a realizability interpretation in $({\mathbf{\Sigma}}^1_2\mbox{-}\mathbf{AC})+\mathbf{BI}$ follows by an extension of the techniques used in \cite[Theorem 6.2]{lpo}.
The proof furnished
a realizability model for $\CZF+\RDC+\LPO$ that is based on recursion in the type-2 object
$E:(\N\to \N)\to \N$ with $E(f)=n+1$ if $f(n)=0$ and $\forall i<n\,f(n)>0$ and $E(f)=0$ if $\forall n\,f(n)>0$.
Recursion in $E$ is
 formalizable in the theory of bar induction, i.e. $(\mathbf{\Pi}^0_{\infty}\mbox{-}\mathbf{CA})+\mathbf{BI}$, which is known to have the same strength as $\CZF$ (see \cite[Theorem 2.2]{lpo}).
 The same recursion theory (or partial combinatory algebra) can be employed in extending the modeling of  a type structure
 given in \cite[$\S5$]{lpo} to the larger type structure needed for $\CZF+\REA+\RDC+\LPO$. This is achieved by basically taking the type structure in \cite[5.8]{rohio} but changing the underlying partial combinatory algebra to the one obtained from recursion in the type two object $E$ rather than the usual one provided by the partial recursive functions on $\mathbb{N}$.

 It is very likely that the interpretation also validates $\mathbf{\Pi\Sigma}{\mathsf{W}}\mbox{-}\AC$ and $\PAx$, but this hasn't yet been checked.

 At any rate, we have shown the proof-theoretic equivalence of all theories.
\qed
The foregoing proof establishes the claimed results, however, we'd like to look  at Theorem \ref{Theorem1} in more detail, especially at its proof(s)  and the information one can extract from it.

 For starters, what does the phrase ``same proof-theoretic strength" mean?
At a minimum  it means that the theories ought to be finitistically equiconsistent. Here it means that they prove at least the same $\Pi^0_2$ statements of the language of first-order arithmetic. But more can be shown.
A result we will be working toward is that many of the intuitionistic theories of
  Theorem \ref{Theorem1} prove the same arithmetical statements. In particular it will be shown that the extensional and the intensional type theories prove the same arithmetical statements.
  An arithmetical statement gives rise to a type via the propositions-as-types paradigm, so by conservativity of one type theory over another
  with respect to arithmetic statements we mean that the same arithmetical types are provably inhabited in
  both theories.

  The question of the relation between intensional and extensional type theories has been addressed before by Hofmann
  in \cite{hofmann}. The set-up there, though, is somewhat different in that the intensional type theory $TT_I$ of \cite{hofmann} is not a pure intensional type theory. It
  has two extensional rules called {\em functional extensionality} and {\em uniqueness of identity}:
 \begin{eqnarray*} (\mbox{ID-UNI-I})&&\quad{\Gamma\vdash s\intt A\qquad \Gamma\vdash p\intt \bI(A,s,s)
\over \Gamma\vdash IdUni(A,s,p)\intt \bI(\bI(A,s,s),p,\refl(s))}
\\ &&
\vphantom{ffffffffff}\\
(\mbox{EXT-FORM})&&\quad{\Gamma\vdash f,g\intt \Pi_{(x\intt A)}B(x)\qquad \Gamma,\,x\intt A\vdash p(x)\intt \bI(B(x),fx,gx)
\over \Gamma\vdash Ext(f,g,p)\intt \bI(\Pi_{(x\intt A)}B(x),f,g)\,.}
\end{eqnarray*}
These rules are not provable in the purely intensional context, so
as a result, we are pursuing a different question here.

\begin{prop}\label{prop1} $\TO^i$ can be interpreted in $\CZF+\REA$. The interpretation preserves (at least) all arithmetic statements.\end{prop}
\prf The proof of \cite{rohio} Theorem 3.9 provides an interpretation of $\TO^i$ in $\CZF+\REA$
which is essentially a class model of $\TO^i$ inside  $\CZF+\REA$. Having defined an applicative structure, the classifications are defined inductively along the (intuitionistic) ordinals. This is inspired by Feferman's construction of a model of $\TO^i$ in \cite[Theorem 4.1.1]{Fef75}. Inspection of the translation confirms that arithmetic statements get preserved. \qed

\begin{prop}\label{prop2}
\begin{itemize}
\item[(i)] $\CZF +\REA$ has an interpretation in  $\mathbf{MLT}_{1\mathsf{W}}\mathsf{V}$.
\item[(ii)] $\CZF+\REA+{\mathbf{\Pi\Sigma}}{\mathsf{W}}\mbox{-}\AC +\RDC + \PAxx$ has an interpretation in
$\mathbf{MLT}_{1\mathsf{W}}^{\ext}\mathsf{V}$.
 \end{itemize}
\end{prop}
\prf (i) and (ii) follow from \cite{aczel86}. The interpretation uses the type $\mathsf{V}$ and two propositional
functions \begin{eqnarray*} \dot{=} \,:\,\mathsf{V}\times \mathsf{V}\to \mathcal{U}_0 && \\
\dot{\in} \,:\,\mathsf{V}\times \mathsf{V}\to \mathcal{U}_0
\end{eqnarray*} to interpret $=$ and $\in$.
For (i), the identity type does not play any role. For (ii) one needs the extensionality of function types.
\qed

\begin{prop}\label{prop2a}
 $\CZF+\REA+{\mathbf{\Pi\Sigma}}{\mathsf{W}}\mbox{-}\AC +\RDC + \PAxx$ is conservative over $\CZF +\REA+\FTAC$
 for statements of finite type arithmetic (i.e., of the language of $\mathbf{HA}^{\omega}$).
\end{prop}
\prf From \cite[Theorem 5.23]{rt} it follows that $\CZF+\REA+{\mathbf{\Pi\Sigma}}{\mathsf{W}}\mbox{-}\AC +\RDC + \PAxx$
and $\CZF+\REA+{\mathbf{\Pi\Sigma}}{\mathsf{W}}\mbox{-}\AC$ prove the same sentences of finite type arithmetic (and more)
since the inner model $\mathsf{H}(\mathsf{Y}^*_{\mathsf{W}})$ satisfies $\CZF+\REA+{\mathbf{\Pi\Sigma}}{\mathsf{W}}\mbox{-}\AC +\RDC + \PAxx$, assuming $\CZF+\REA+{\mathbf{\Pi\Sigma}}{\mathsf{W}}\mbox{-}\AC$ in the background.

By \cite[Theorem 4.33]{ober}, there is an interpretation of $\CZF+\REA+{\mathbf{\Pi\Sigma}}{\mathsf{W}}\mbox{-}\AC$ in $\CZF+\REA$.
Inspection shows that, in the presence of $\FTAC$, the meanings of statements of finite type arithmetic are
 preserved under this interpretation. \qed

\begin{prop}\label{prop3} For $\theta$ a sentence of arithmetic let $\|\theta\|$ be the corresponding  type term
according to the propositions-as-types translation.
If $$\mathbf{MLT}_{1\mathsf{W}}^{\ext}\mathsf{V}\vdash t\,:\,\|\theta\|$$ for some term $t$, then
$$\CZF+\REA+ \FTAC\vdash \theta_{\mbox \tiny set}$$ with $\theta_{\mbox \tiny set}$ denoting the standard set-theoretic rendering of $\theta$.
\end{prop}
\prf Assume $\mathbf{MLT}_{1\mathsf{W}}^{\ext}\mathsf{V}\vdash t\,:\,\|\theta\|$. The interpretation $^\wedge$ of $\mathbf{MLT}_{1\mathsf{W}}^{\ext}\mathsf{V}$ into $\CZF+\REA$ given in \cite[$\S6$]{rt} yields
$\CZF+\REA\vdash (t\,:\,\|\theta\|)^\wedge$. Inspection shows that $(t\,:\,\|\theta\|)^\wedge$ is a statement about the finite type structure over $\omega$. One then sees,  with the help of $\FTAC$, that
$\theta_{\mbox \tiny set}$ holds. This is similar to the proof of \cite[Theorem 3.15]{rt}.
 \qed

 \begin{thm}\label{ikpi}
  $\CZF+\REA+{\mathbf{\Pi\Sigma}}{\mathsf{W}}\mbox{-}\AC +\RDC + \PAxx$ is conservative over $$\mathbf{IKP}+\forall x\,\exists y\,[x\in y\;\wedge\;\mbox{\tt $y$ is an admissible set}]$$
 for arithmetical statements.
 \end{thm}

 \prf We shall use the shorthand $\mathbf{IKPi}$ for the latter theory. By Proposition \ref{prop2a} it suffices to show that $\CZF +\REA+\FTAC$ is conservative over $\mathbf{IKPi}$ for arithmetic statements.
\cite[Theorem 5.11]{rohio} shows that $\mathbf{MLT}_{1\mathsf{W}}^{\ext}\mathsf{V}$ has an interpretation in the classical theory $\mathbf{KPi}$ where
types are interpreted as subsets of $\omega$ and crucially dependent products of types   are interpreted as sets of indices of partial recursive functions. This also furnishes an interpretation of $\CZF+\REA+\FTAC$ in $\mathbf{KPi}$ since the former is interpretable in $\mathbf{MLT}_{1\mathsf{W}}^{\ext}\mathsf{V}$. The interpretation also works for $\mathbf{IKPi}$ as definition by (transfinite) $\Sigma$-recursion
works in intuitionistic $\mathbf{KP}$ as well (see \cite[Sec. 11]{maralt} and \cite[Sec. 19]{mar}). The inductive
definition of 5.8 in \cite{rohio} proceeds along the ordinals and focusses on successor ordinals, seemingly requiring
 a classical case distinction as to whether an ordinal is a successor or a limit or 0, but this is actually completely irrelevant.

 Now, the upshot of this hereditarily recursive interpretation is that every $\Pi^0_2$ theorem of $\CZF+\REA+\FTAC$ is provable in $\mathbf{IKPi}$. To be able to extend this approach to all of arithmetic, one needs a more abstract type structure such that interpretability entails deducibility. The conservativity of $\mathbf{HA}^{\omega}+\FTAC$ over $\mathbf{HA}$, due to Goodman
  \cite{goodman,Goodman}, provides the template.  The two steps of Goodman's second proof have been neatly separated by Beeson \cite{Beeson} to construct a general methodology for showing an intuitionistic theory $T$ to be conservative over another theory $S$ for arithmetic statements. The idea is
to combine two  interpretations, where the first uses functions that are recursive relative to a generic oracle and
the second step is a forcing construction. The same idea has been used by Gordeev \cite{gordeev}, and in more recent times by Chen and Rathjen in \cite{rmc,rr,turingR},   establishing several conservativity results.

 The oracle $\mathcal O$ will be a fixed but arbitrary partial function
from $\mathbb N$ to $\{0,1\}$.  A partial function $\phi$ is recursive relative to $\mathcal O$ if it is given by a Turing machine with
 access to $\mathcal O$. During a computation the oracle  may be consulted about the value of $\mathcal O(n)$ for several $n$. If $\mathcal O(n)$ is defined it will return that value and the computation will continue, but if $\mathcal O(n)$ is not defined no response
will be coming forward and the computation will never come to a halt. The idea of the second interpretation step is that on account of $\mathcal O$'s arbitrariness it can be interpreted in many ways. Given
an arithmetic statement $\theta$, an oracle ${\mathcal O}_\theta$ can be engineered so that in a forcing model realizability of $\theta$
with functions computable relative to ${\mathcal O}_\theta$ entails
the truth of $\theta$. The final step, then,  is achieved by noticing that for arithmetic statements forcibility (where the forcing conditions are finite partial functions on $\mathbb N$)
 and validity coincide.
 For details we'll have to refer to \cite{rmc,rr}. \qed

 \begin{deff}{\em Below we shall speak about arithmetical statements in various theories with differing languages.
 There is a canonical translation of the language of first and second order arithmetic into the language of set theory.
 However, it is perhaps less obvious what arithmetical statements mean in the context of type theory.

 The terms of the language of $\HA$ are to be translated in an obvious way, crucially using the type-theoretic recursor for the
 type $\mathbb{N}$. In this way each term $t$ of $\HA$ gets assigned a raw term $\hat{t}$ of type theory. For details see \cite[pp. 71--75]{ml84}, \cite[XI.17]{beeson} \cite[Ch. 11, Sect. 2]{TD88}.
 An equation $s=t$ of the language $\HA$ is translated as a type-expression $\mathsf{Id}(\mathbb{N},\hat{s},\hat{t})$.
 For complex formulas the translation proceeds in the obvious way.

 We then say that two type theories $TT_1$ and $TT_2$ prove the same arithmetical statements if for all sentences $A$ of $\HA$,
  \begin{eqnarray*} TT_1\vdash p\,:\,\hat{A}\mbox{ for some $p$ }&\mbox{ iff }& TT_2\vdash p'\,:\,\hat{A}\mbox{ for some $p'$, }
 \end{eqnarray*} where $\hat{A}$ denotes the type-theoretic translation of $A$.}\end{deff}

Recall that  $\mathbf{IKPi}$ is the theory $\mathbf{IKP}+\forall x\,\exists y\,[x\in y\;\wedge\;\mbox{\tt $y$ is an admissible set}].$

\begin{thm}\label{Theorem2} The following theories prove the same arithmetical statements, i.e. statements of the language of first order arithmetic (also known as Peano arithmetic).
\begin{itemize} \item[(i)]  $\TO^i$.
\item[(ii)]  $\CZF+\REA$.
\item[(iii)] $\CZF+\REA+{\mathbf{\Pi\Sigma}}{\mathsf{W}}\mbox{-}\AC +\RDC + \PAxx$.
\item[(iv)]  $\mathbf{MLT}_{1\mathsf{W}}^{\ext}\mathsf{V}$.
 \item[(v)]   $\mathbf{MLT}_{1\mathsf{W}}\mathsf{V}$.
\item[(vi)]  $\mathbf{MLT}_{1\mathsf{W}}^{\ext}$.
\item[(vii)] $\mathbf{MLT}_{1\mathsf{W}}$.
\item[(viii)] $\mathbf{IARI}$.
\item[(ix)] $\mathbf{IKPi}$.
\end{itemize}
\end{thm}
\prf Let $\theta$ be an arithmetic sentence. Then we have
\begin{eqnarray*}\TO^i\vdash \theta &\Rightarrow & \CZF+\REA\vdash\theta \\ &
\Rightarrow& \CZF+\REA+{\mathbf{\Pi\Sigma}}{\mathsf{W}}\mbox{-}\AC +\RDC + \PAxx\vdash\theta\\
&\Rightarrow & \mathbf{IKPi}\vdash \theta\end{eqnarray*}
by Proposition \ref{prop1} and Theorem \ref{ikpi}.
Now it follows from J\"ager's article \cite{j83} and from \cite{jp} that every initial segment of the proof-theoretic ordinal of $\mathbf{IKPi}$
is provably well-founded in $\TO^i$, and thus, if $\mathbf{IKPi}\vdash \theta$, then $\TO^i$ is sufficient to show that
there is an infinite intuitionistic cut-free proof of $\theta$. By induction on the length of the proof it then follows that all sequents in the proof are true, yielding that $\TO^i\vdash \theta$.
The upshot is that the theories of (i), (ii) and (iii) prove the same arithmetic statements. Furthermore, if
$\mathbf{MLT}_{1\mathsf{W}}^{\ext}\mathsf{V}\vdash t\,:\,\|\theta\|$ for some term $t$, then $\CZF+\REA\vdash \theta$ by Proposition \ref{prop3} and hence
$\TO^i\vdash \theta$.

So to finish the proof it would suffice to show that
 $\TO^i\vdash \theta$ yields $\mathbf{MLT}_{1\mathsf{W}}\mathsf{V}\vdash s\,:\,\|\theta\|$ for some term $s$.
 Now \cite[Sec. 6]{rohio} shows that the intuitionistic  theory $\mathbf{IARI}$ has the same proof-theoretic ordinal as $\mathbf{IKPi}$ and $\TO^i$. So from $\TO^i\vdash \theta$ it follows that
 $\mathbf{IARI}\vdash \theta$. By \cite[Theorem 6.9]{rohio} we then get $\mathbf{MLT}_{1\mathsf{W}}\mathsf{V}\vdash s\,:\,\|\theta\|$ for some term $s$, completing the circle.
 \qed

 \begin{remark}{\em Ordinal analysis played a crucial role in the proofs of Theorems \ref{Theorem1} and \ref{Theorem2}.
 Having the same proof-theoretic ordinal allowed us to infer that $\TO^i$, $\mathbf{IKPi}$ and $\mathbf{IARI}$ prove the same arithmetic statements.

  For a long time \cite{j83} was also the only proof that enabled one to reduce the classical theories
  $({\mathbf{\Delta}}^1_2\mbox{-}\mathsf{CA})+\mathsf{BI}$ and $\mathbf{KPi}$ to classical $\TO$. There is now also a proof by Sato \cite{sato} for the reductions in the {\em classical} case that avoids proof-theoretic ordinals.  However, determining the strength of other important fragments of $\mathbf{MLTT}$ (such as the ones analyzed by Setzer in \cite{set}) still requires the techniques of 
ordinal analysis.

The strength of other important fragments of $\mathbf{MLTT}$ was analyzed by Setzer in \cite{set}. }
  \end{remark}

  \begin{remark}{\em We conjecture that also the theory $\CZF+\sREA+{\mathbf{\Pi\Sigma}}{\mathsf{W}}\mbox{-}\AC +\RDC + \PAxx$
  (or at least $\CZF+\sREA+{\mathbf{\Pi\Sigma}}{\mathsf{W}}\mbox{-}\AC +\RDC$) proves the same arithmetical statements as any of the theories featuring in Theorem \ref{Theorem2}. As the latter relies on a substantial number of results from the literature, several of them would have to be revisited and possibly amended to establish this.}
  \end{remark}

\section{On relating theories II: $\MLTT$ and friends}
So far we have only gathered results concerning theories that are of the strength of Martin-L\"of type theory with one universe. The earlier quote by Harris speculated on the strength of type theory with infinitely many universes.
As it turns out, similar techniques can be applied in this context as well.

To begin with, we shall define versions of explicit mathematics, second order arithmetic and constructive set theory featuring analogues of universes.
\subsection{$\TO^i$ with universes.}
\begin{deff}{\em Systems of explicit mathematics with universes have been defined and
studied in several papers (cf. \cite{js,jstu,jks}) and were probably first introduced by Feferman \cite{fef-hancock}.

By $\TO^i+\bigcup_n\mathsf{U}_n$ we denote an extension of $\TO^i$ whose language has infinitely many classification constants $\mathsf{U}_0,\mathsf{U}_1,\ldots$ and the following axioms for
each constant $\mathsf{U}_n$.
\begin{enumerate}
\item $\mathbf N\vin \mathsf{U}_n$ and $\mathsf{U}_i\vin \mathsf{U}_n$ for $i<n$.
\item  $\forall x\vin \mathsf{U}_n\, \exists X\, x=X$ (i.e. every element of $\mathsf{U}_n$ is a classification).
    \item For every elementary formula  $\psi(x,\vec v, X_1,\ldots,X_r)$  with all  classification variables exhibited and which does not contain constants $\mathsf{U}_i$ with $i\geq n$, $$\forall X_1,\ldots,X_r\vin \mathsf{U}_n\, \exists Y\,[Y\vin \mathsf{U}_n\,\wedge\,
    Y\simeq\{x\,:\,\psi(x,\vec v, X_1,\ldots,X_r)\}]\,.$$
    \item $\forall X\vin \mathsf{U}_n\,[\forall x\vin X\,
\exists Y\vin \mathsf{U}_n\ fx\simeq Y \rightarrow \exists Z[Z\in \mathsf{U}_n\,\wedge\,Z \simeq \mathbf{j}(X,f)]]$.

\item $\forall X,Y\vin \mathsf{U}_n\, \exists Z\,[Z\in \mathsf{U}_n\,\wedge\,Z \simeq \bi(X,Y)]$.
\end{enumerate}
In other words, a classification $\mathsf{U}_n$ is a universe containing $\mathbf{N},\mathsf{U}_0,\ldots, \mathsf{U}_{n-1}$ closed under elementary comprehension, join and inductive generation.\\[1ex]
By $\TO^i+\bigcup_{i<n}\mathsf{U}_{i}$ we denote the theory with just the universes $\mathsf{U}_0,\ldots,\mathsf{U}_{n-1}$ and
their pertaining axioms.
}\end{deff}

\subsection{Universes in intuitionistic second order arithmetic.}
It is also useful to have a many universes version of \Iari to obtain an intuitionistic theory of second order arithmetic which can be easily interpreted in $\MLTT$. One idea would be to adopt the notion of $\beta$-model from Definition \ref{4.8} to serve as a notion of universe. However, a $\beta$-model comes with an explicit countable enumeration of its sets and therefore it would be difficult if not impossible
to model such structures in $\MLTT$.
Instead, an option is to add set predicates ${\mathfrak U}_0, {\mathfrak U}_1, \ldots$ to the language ${\mathcal L}_2$ that are intended to apply to sets of natural numbers with the aim of
 singling out collections of sets that have universe-like properties.

 \begin{deff}{\em
 The theory \Iario+$\bigcup_n{\mathfrak U}_n$ has additional predicates  ${\mathfrak U}_0, {\mathfrak U}_1, \ldots$ for
 creating new atomic formulas ${\mathfrak U}_n(X)$ ($n\in\mathbb{N}$),
 where $X$ is a second order variable. We use abbreviations like $\forall X\in {\mathfrak U}_n\,\varphi$ and $\exists X\in {\mathfrak U}_n\,\varphi$ for $\forall X\,({\mathfrak U}_n(X) \to\varphi)$ and $\exists X\,({\mathfrak U}_n(X)\,\wedge\,\varphi)$, respectively.
 If $\psi$ is any formula of this language, then $\psi^{{\mathfrak U}_n}$ arises from $\psi$
 by relativizing all second order quantifiers to ${\mathfrak U}_n$, i.e., replacing all quantifiers  $QX$ in $\psi$ by $QX\in {\mathfrak U}_n$.

   In addition to the axioms of \Iari there are the following pertaining to the new predicates.

\begin{enumerate}
\item The predicates ${\mathfrak U}_n$ are cumulative, i.e.
$\forall X\,[{\mathfrak U}_i(X) \to {\mathfrak U}_j(X)]$ whenever $i\leq j$.
\item {\bf Induction}:
$$\phi(0)\land\forall u[\phi(u)\then\phi(u+1)]\then\forall u\phi(u)$$
for all formulae $\phi$.
\item {\bf Arithmetic Comprehension Schema for ${\mathfrak U}_n$}:
$$Y_1,\ldots,Y_r\in{\mathfrak U}_n\to\exists X\in{\mathfrak U}_n\forall u[u\in X\lra\psi(u,Y_1,\ldots,Y_r)]$$
if $\psi(u,Y_1,\ldots,Y_r)$ is a formula with all free second order variables exhibited, in which all second order quantifiers are of the form $QX\in {\mathfrak U}_i$ for some $i<n$, and moreover, no predicates ${\mathfrak U}_j$ for $j\geq n$ occur in it.

\item {\bf Replacement}: %Axiom of Countable Collection:
$$\forall X\in{\mathfrak U}_n[\forall u\in X\exists\,!
Y\in{\mathfrak U}_n\phi(u,Y)\then\exists Z\in{\mathfrak U}_n\forall u\in X\,\phi(u,
 Z_u)]$$
for all formulas $\phi$. Here $\phi(u,Z_u)$ arises from $\phi(u,Z)$ by
replacing each occurrence $t\in Z$ in the formula by $\pair ut\in Z$.
\item {\bf Inductive Generation}:
$$\forall U\in{\mathfrak U}_n\forall X\in {\mathfrak U}_n\exists Y\in {\mathfrak U}_n\,\bigl[
\WPU(X,Y)\land(\forall u[\forall v(v\ltx u\then\phi(v))\then\phi(u)]\then
\forall x\in Y\,\phi(x))\bigr],$$
for all formulas $\phi$, where $v\ltx u$ abbreviates $\pair vu\in X$ and
$\WPU(X,Y)$ stands for
$$\progU(X,Y)\land\forall Z[\progU(X,Z)\then Y\subseteq Z]$$
with $\progU(X,Y)$ being $\forall y \in U[\forall z(z\ltx y\then z\in Y)
\then y\in Y]$.
\end{enumerate}
By \Iari$+\bigcup_{i<m}\mathfrak{U}_{i}$ we denote the theory with only the additional predicates ${\mathfrak U}_0,\ldots,{\mathfrak U}_{m-1}$  and their pertaining axioms.
%axioms the above axioms pertaining to the predicates ${\mathfrak U}_0,\ldots,{\mathfrak U}_{n-1}$.
}\end{deff}

\begin{deff}{\em Recall the notion of inaccessible set defined in \ref{4.7}. For $n>0$, $\mathbf{Inacc}(n)$ stands for the set-theoretic statement that there are $n$-many inaccessible sets $I_0\in\ldots \in I_{n-1}$.
$\mathbf{Inacc}(0)$ stands for $0=0$. 

$\beta$-models were introduced in \ref{4.8}. By $\mathbf{Beta}(n)$ we denote the statement of second order arithmetic asserting that there are $n$ many sets
$A_0,\ldots,A_{n-1}$ which are $\beta$-models of ${\mathbf{\Sigma}}^1_2\mbox{-}\mathbf{AC}$ such that
$A_0\in \ldots \in A_{n-1}$, where for sets $X,Y$ of natural numbers $X\in Y$ is defined by $\exists u\;X=Y_u$.}
\end{deff}

 For $n>0$, let $\mathbf{MLT}_{n\mathsf{W}}V$ be the fragment of $\MLTT$ with $n$-many universes $\mathcal{U}_0,\ldots,\mathcal{U}_{n-1}$, where the $\mathsf{W}$-constructor can solely be applied to families of types in $\mathcal{U}_0,\ldots,\mathcal{U}_{n-1}$ but one can also form  the type $\mathsf{V}:=\mathsf{W}_{(A:\mathcal{U}_{n-1})}A$, i.e. a $\mathsf{W}$-type over the largest universe $\mathcal{U}_{n-1}$.  We shall also consider the type theory $\mathbf{MLT}_{n\mathsf{W}}$ which is the fragment of $\mathbf{MLT}_{n\mathsf{W}}V$ without the type $\mathsf{V}$.
 
 Below we assume that $n>0$.

\begin{thm}\label{vor}
\begin{itemize}
\item[(i)] $\TO^i+\bigcup_{i<n}\mathsf{U}_{i}$ has an interpretation in $\CZF+\REA+\mathbf{Inacc}(n)$. The interpretation preserves (at least) all arithmetic statements.
    \item[(ii)] $\CZF +\REA+\mathbf{Inacc}(n-1)$ has an interpretation in  $\mathbf{MLT}_{n\mathsf{W}}\mathsf{V}$.
\item[(iii)] $\CZF+\REA+\mathbf{Inacc}(n-1)+{\mathbf{\Pi\Sigma}}{\mathsf{W}}\mbox{-}\AC +\RDC + \PAxx$ has an interpretation in
$\mathbf{MLT}_{n\mathsf{W}}^{\ext}\mathsf{V}$.

\item[(iv)] $\mathbf{MLT}_{n\mathsf{W}}^{\ext}\mathsf{V}$ has an interpretation in the classical set theory $\mathbf{KPi}$ plus an axiom asserting that there exist $n-1$-many recursively inaccessible ordinals.

 \item[(v)]    $({\mathbf{\Sigma}}^1_2\mbox{-}\mathbf{AC})+\mathbf{BI}+\mathbf{Beta}(n)$ has
 an interpretation in $\mathbf{KPi}$ plus the existence of $n$-many recursively inaccessible ordinals.

 \item[(vi)] $\mathbf{KPi}$ plus the existence of $n$-many recursively inaccessible ordinals has a sets-as-trees interpretation in $({\mathbf{\Sigma}}^1_2\mbox{-}\mathbf{AC})+\mathbf{BI}+\mathbf{Beta}(n)$.

 \item[(vii)]   The intuitionistic system $\mathbf{IRA}+\bigcup_{i<n-1}\mathfrak{U}_{i}$ of second order arithmetic can be interpreted
      in $\mathbf{MLT}_{n\mathsf{W}}$.
  \item[(viii)]    $\CZF+\REA+\RDC + \mathbf{Inacc}(n)$ has a realizability interpretation in
  $\mathbf{KPi}$ plus the existence of $n$-many recursively inaccessible ordinals.
 \item[(ix)] All the above theories have the same proof-theoretic strength and prove (at least) the same $\Pi^0_2$-statements of arithmetic.

     \end{itemize}

\end{thm}
\prf The interpretations are extensions of those discussed in the previous section, taking more universes into account.
We can only indicate the steps. The interpretation of $\TO^i$ in $\CZF+\REA$ can be lifted to an interpretation of $\TO^i+\bigcup_{i< n}\mathsf{U}_i$  into   $\CZF+\REA+\mathbf{Inacc}(n)$. The latter theory possesses a sets-as-types interpretation in intensional
Martin-L\"of type theory with $n+1$ universes.

$\CZF+\REA+{\mathbf{\Pi\Sigma}}{\mathsf{W}}\mbox{-}\AC +\RDC + \PAxx+\mathbf{Inacc}(n-1)$ possesses a sets-as-types interpretation in $\mathbf{MLT}_{n\mathsf{W}}^{\ext}\mathsf{V}$. In turn, $\mathbf{MLT}_{n\mathsf{W}}^{\ext}\mathsf{V}$
 can be interpreted in  classical Kripke-Platek set theory $\mathbf{KPi}$ plus an axiom asserting that there are at least  $n-1$-many recursively inaccessible ordinals, following the Ansatz of \cite[Theorem 5.11]{rohio}. $({\mathbf{\Sigma}}^1_2\mbox{-}\mathbf{AC})+\mathbf{BI}+\mathbf{Beta}(n)$
can be easily interpreted in $\mathbf{KPi}$ plus $n$-recursively inaccessible ordinals.

%A further crucial ingredient  is an ordinal analysis of Kripke-Platek set theory $\mathbf{KP}$ plus for every $n>0$ an axiom %asserting that there are at least  $n$-many recursively inaccessible ordinals.
%It can be shown that for fixed $n$ the ordinal analysis of $\mathbf{KP}$ plus $n$-many
%recursively inaccessible ordinals can be carried out in both $\TO^i+\bigcup_{i\leq n+1}\mathsf{U}_i$
%and \Iario$+\bigcup_{i\leq n+1}{\mathfrak U}_i$.

The proof-theoretic equivalence ensues from an ordinal analysis of the `top theory',
$\mathbf{KPi}$ plus the existence of $n$-many recursively inaccessible ordinals, together with proofs that
any ordinal below the proof-theoretic ordinal of that theory is provably well-founded in
$\TO^i+\bigcup_{i<n}\mathsf{U}_{i}$ as well as \Iari$+\bigcup_{i<n}\mathfrak{U}_{i}$.
Neither the ordinal analysis nor the well-ordering proofs are available from the published literature. The ordinal analysis of $\mathbf{KPi}$ plus the existence of $n$-many recursively inaccessible ordinals, though, can be obtained in a straightforward way by extending the one given for $\KPi$ in \cite{jp} or rather its modern version in \cite{bu93}. It also follows from the ordinal analysis of the much stronger theory $\mathbf{KPM}$ given in \cite{r91} by restricting the treatment therein to
the pertaining small fragments. For the well-ordering proof substantially more work is required; details will be published in \cite{rth}.
\qed

\begin{thm}\label{Theorem3} The following theories have the same proof-theoretic strength and prove the same $\Pi^0_2$-statements of arithmetic:
\begin{itemize}
\item[(i)] $\TO^i+\bigcup_n\mathsf{U}_n$.
\item[(ii)]  $\CZF$ plus $\mathbf{Inacc}(n)$ for all $n>0$.
\item[(iii)] $\CZF+{\mathbf{\Pi\Sigma}}{\mathsf{W}}\mbox{-}\AC +\RDC + \PAxx$ plus
$\mathbf{Inacc}(n)$ for all $n>0$.
\item[(iv)] The extensional type theory $\mathbf{MLTT}^{\ext}$.
 \item[(v)]   $\mathbf{MLTT}$.
\item[(vi)] The classical subsystem of second order arithmetic $({\mathbf{\Sigma}}^1_2\mbox{-}\mathbf{AC})+\mathbf{BI}$
plus $\mathbf{Beta}(n)$ for all $n>0$.
\item[(vii)] Classical Kripke-Platek set theory $\mathbf{KP}$ plus for every $n>0$ an axiom asserting that there are at least  $n$-many recursively inaccessible ordinals.
 \item[(viii)]    \Iario$+\bigcup_n{\mathfrak U}_n$.
\item[(ix)] $\CZF+\RDC + \LPO$ plus the axioms $\mathbf{Inacc}(n)$ for all $n>0$.
\end{itemize}

%\Pi\Sigma{\mathsf{W}}\mathsf{I}_0\ldots\mathsf{I}_n\mbox{-}\AC
\end{thm}

\prf This follows directly from the previous theorem. \qed

The latter theorem also shows that the strength of $\MLTT$ is dwarfed by that of $(\Pi^1_2\mbox{-}\mathbf{CA})$.
It corresponds to a tiny fragment of second order arithmetic which itself is a tiny fragment of $\ZF$, so there are aeons between $\MLTT$ and
 classical set theory with inaccessible cardinals.

\begin{thm}\label{Theorem4} The following theories prove the same arithmetical statements:
\begin{itemize}
\item[(i)] $\TO^i+\bigcup_n\mathsf{U}_n$.
\item[(ii)]   $\mathbf{MLTT}$.
\item[(iii)] The extensional type theory $\mathbf{MLTT}^{\ext}$.
\item[(iv)]  $\CZF$ plus $\Inac(n)$ for every $n>0$.
\item[(v)] $\CZF+\bigcup_n\Inac(n)+\bigcup_n{\mathbf{\Pi\Sigma}}{\mathsf{W}}\mbox{-}\AC +\RDC + \PAxx$.
\item[(vi)] \Iario$+\bigcup_n{\mathfrak U}_n$.

 \end{itemize}
\end{thm}
\prf The methods for proving this were described in the proof of \ref{Theorem2}. Details will appear in \cite{rth}. \qed

Finally, it should be mentioned that Martin-L\"of type theory with stronger universes (e.g. Mahlo universes) has been studied by Setzer (cf. \cite{setz}).

\subsection{Adding the Univalence Axiom}
The quote (\ref{hybris}) from Harris' book \cite{harris} claimed that modeling  Voevodsky's univalence axiom ($\UA$) requires infinitely many inaccessible cardinals (for a definition of $\UA$ see \cite[Sec. 2.10]{hott}).
While the simplicial model of type theory with univalence developed in the paper \cite{klv} by Kalpulkin, Lumsdaine and Voevodsky is indeed carried out in a background set theory with inaccessible cardinals, it is by no means clear that the existence or proof-theoretic strength of these objects is required for finding a model of type theory with $\UA$. In actuality, Bezem, Coquand and Huber in their article \cite{bch} provided a cubical model of type theory that also validates $\UA$. Crucially,
their modeling can be carried out in a constructive background theory such as $\CZF+\bigcup_n\Inac(n)+\bigcup_n{\mathbf{\Pi\Sigma}}{\mathsf{W}}\mbox{-}\AC +\RDC + \PAxx$.
Thus it follows that adding $\UA$ does not increase the strength of type theory and that no inaccessible cardinals are required. Hence in view of Theorem \ref{Theorem3}  we have the following result.

\begin{cor}\label{UA1} $\mathbf{MLTT}$ has the same proof-theoretic strength as $\mathbf{MLTT}+\UA$.
Thus $\mathbf{MLTT}+\UA$ shares the same proof-theoretic strength with all theories listed in Theorem \ref{Theorem3}, in particular with
classical Kripke-Platek set theory $\mathbf{KP}$ augmented by axioms  asserting that there are at least  $n$-many recursively inaccessible ordinals for every $n>0$.
\end{cor}

\section{On relating theories III: Omitting $\mathsf{W}$}
The proof-theoretic strength of type theories crucially depends on the availability of inductive types and to a much lesser extent on its
universes. Relinquishing the $\mathsf{W}$-type brings about an enormous collapse of proof power (cf. \cite{super,saper,r97b}). Letting $\mathbf{MLTT}^-$ be $\mathbf{MLTT}$ bereft of the $\mathsf{W}$-type constructor, we arrive at a theory no stronger than the system $\mathbf{ATR}_0$ of reverse mathematics (see \cite[I.11]{simpson}), having the famous ordinal $\Gamma_0$
as its proof-theoretic ordinal. According to Feferman's analysis (see \cite{fe64,f68}), $\Gamma_0$ delineates the limit of a notion of predicativity
that only accepts the natural numbers as a completed infinity (which was first adumbrated in Hermann Weyl's book {\em ``Das Kontinuum"} from 1918 \cite{weyl}).    Peter Hancock conjectured in the 1970s  the ordinal of $\mathbf{MLTT}^-$
to be $\Gamma_0$.
 Feferman \cite{fef-hancock} and independently Aczel (see also \cite{aa})  proved
  {\em Hancock's Conjecture}.
  There is also a version of $\CZF$ with inaccessible sets of strength $\Gamma_0$,  due to Crosilla and Rathjen \cite{crosilla-rathjen},  which does not have set induction. Thus the set-theoretic analogue to eschewing $\mathsf{W}$-types consists in leaving out the principle of set induction.
  In the next theorem we denote by $\mathbf{ATR}_0^i$ the intuitionistic version of $\mathbf{ATR}_0$ (see \cite[Definition 4.10]{super} for details). By $\CZF^-$ we denote Constructive Zermelo-Fraenkel set theory without set induction but with the Infinity axiom strengthened as follows:
  \begin{eqnarray}\label{inf1} && 0\in \omega \,\wedge\,\forall y[\,y\in \omega\to y+1\in \omega] \\ \label{inf2}
  && \forall x\,[\,0\in x\,\wedge\,\forall y(y\in x\to y+1\in x) \to \omega\subseteq x]
  \end{eqnarray} (for details see \cite[Definition 2.2]{crosilla-rathjen}).
  Likewise we denote by $\KP^-$ the theory without the set induction scheme but with the infinity axioms
  (\ref{inf1}) and (\ref{inf2}).

  The notion of weak inaccessibility used below is the one from Definition \ref{4.4}.
  For $n>0$ let $\mathbf{wInacc}(n)$ be the statement that there exist weakly inaccessible  sets $x_0,\ldots,x_{n-1}$ such that $x_0\in \ldots \in x_{n-1}$.

  A restricted form of $\RDC$ is $\Delta_0\mbox{-}\RDC$: For all $\Delta_0$-formulae
$\theta$ and $\psi$, whenever
$$(\forall x\in a)\bigl[\theta(x)\,\rightarrow\,(\exists y\in a)
\bigl(\theta(y)\,\wedge\,\psi(x,y)\bigr)\bigr]$$ and $b_0\in a\,\wedge\,\phi(b_0)$,
then
there exists a function $f:\omega\rightarrow a$ such that $f(0)=b_0$ and
$$(\forall n\in\omega)\bigl[\theta(f(n))\,\wedge\,\psi(f(n),f(n+1))\bigr].$$

  \begin{thm}\label{Theoremminus1}
  The following theories share the same proof-theoretic strength and ordinal $\Gamma_0$, and
  prove the same $\Pi^0_2$-sentences of arithmetic:
  \begin{itemize}
  \item[(i)] $\mathbf{MLTT}^-$.
  \item[(ii)] The extensional version of $\mathbf{MLTT}^-$.
  \item[(iii)]  $\mathbf{ATR}_0$.
  \item[(iv)]  $\mathbf{ATR}_0^i$.
  \item[(v)] $\CZF^-+ \forall x\,\exists y\,[\,x\in y\,\wedge\,\mbox{$y$ is weakly inaccessible}]+\Delta_0\mbox{-}\RDC$.
  \item[(vi)] $\CZF^-+ \{\mathbf{wInacc}(n)\mid n>0\}+\RDC$.
  \item[(vii)] $\KP^-+\forall x\,\exists y\,[x\in y\,\wedge\,\mbox{$y$ is admissible}]$.
  \end{itemize}
  \end{thm}
  \prf We only have to establish  that all theories have proof-theoretic ordinal $\Gamma_0$.
  For extensional $\mathbf{MLTT}^-$ this follows from \cite{fef-hancock}. The lower bound part, namely that  $\mathbf{MLTT}^-$ has at least the strength  $\Gamma_0$ is due to Jervell \cite{jervell}.
  %I think the proofs also work for the intensional theory $\mathbf{MLTT}^-$.
  So we are done with (i) and (ii). That $\mathbf{ATR}_0$ has ordinal $\Gamma_0$ is well known. For $\mathbf{ATR}_0^i$
  this follows from the observation in \cite[Lemma 4.11]{super} that the well ordering proof for any ordinal notation below $\Gamma_0$ uses only intuitionistic logic.
  The determination of the ordinal for the system in (v) and (vi) is due to Crosilla and Rathjen
  \cite[Corollary 9.14]{crosilla-rathjen} with the validation of $\Delta_0\mbox{-}\RDC$ and $\RDC$ coming from
  \cite[Theorem 4.17]{anti} and \cite[Theorem 4.16]{anti}, respectively.
  The proof-theoretic analysis of  the system in (vii) is due to J\"ager \cite{J}.
  \qed

  We also conjecture that all of the intuitionistic theories from the above list, i.e., $\mathbf{MLTT}^-$, the extensional version of $\mathbf{MLTT}^-$, $\mathbf{ATR}_0^i$, and $\CZF^-+ \forall x\exists y\,[ x\in y\,\wedge\, \mbox{ $y$ is weakly inaccessible}]$ prove the same arithmetic statements using the usual techniques. But we have not yet checked that. What is known is that $\mathbf{ATR}_0^i$ embeds in all of these theories (see \cite{super}).

  A final question concerns the status of the univalence axiom. Do we get more strength when we add $\UA$ to $\mathbf{MLTT}^-$? It turns out that we just have to check whether the cubical model construction from \cite{bch} can be carried out in one of the theories from the list. Inspection of \cite{bch} reveals that  $$\CZF^-+ \forall x\exists y\,[ x\in y\,\wedge\, \mbox{ $y$ is weakly inaccessible}]+\Delta_0\mbox{-}\RDC$$ suffices as a background theory for all the constructions, except $\mathsf{W}$-types.

  \begin{cor} The univalent type theory $\mathbf{MLTT}^-+\UA$ is of the same strength as $\mathbf{MLTT}^-$ and $\mathbf{ATR}_0$ and all the other systems from Theorem \ref{Theoremminus1}. Therefore its proof-theoretic ordinal is $\Gamma_0$.
  \end{cor}

\section{Monotone Fixed Point Principles in Intuitionistic Explicit Mathematics}\label{secmon}
Martin-L\"of type theory appears to capture the abstract notion of an inductively defined type very well via its $\mathsf{W}$-type.
There are, however,  intuitionistic theories of inductive definitions that at first glance appear to be just slight extensions
of Feferman's explicit mathematics (see Feferman's quote from Sect. 1) but have turned out to be much stronger than anything considered
in Martin-L\"of type theory. They are obtained from $\TO^i$  by the augmentation of a monotone fixed point
principle which asserts that every monotone operation on
classifications (Feferman's notion of set) possesses a least fixed
point. To be more precise, there are two versions of this principle. $\mathbf{MID}$  merely postulates the
existence of a least solution, whereas $\mathbf{UMID}$ provides a uniform version of this axiom
 by adjoining a new functional  constant to the language,
ensuring that a fixed point is uniformly presentable as a
function of the monotone operation.

\begin{deff}{\em
For {\em extensional equality} of classifications we use the shorthand
``$\dg$'', i.e. $$X\dg Y\ist \forall v(v\elee X\lra v\elee Y).$$ Further, let
$X\extsubseteq Y$ be a shorthand for $\forall v(v\elee X\then v\elee Y).$
To state the monotone
fixed point principle for subclassifications of a given classification $A$
we introduce the following shorthands:
$$\begin{array}{l@{\quad\mbox{if}\quad}l}
\Clop(f,A)&\forall X\,\extsubseteq A\,\,\exists Y\extsubseteq A\,\,
fX\simeq Y\\[2ex]
\Ext(f,A)&\forall X\extsubseteq A\,\,\forall Y\extsubseteq A\,\,[X\dg Y
\then fX\dg fY]\\[2ex]
\Mon(f,A)&\forall X\extsubseteq A\,\,\forall Y\extsubseteq A\,[X\extsubseteq Y
\then fX\extsubseteq fY].\\[2ex]
\Lfp({Y},{f,A})& fY\extsubseteq Y\;\wedge\;Y\extsubseteq A\;\wedge\;
\forall X\extsubseteq A\,\,
\bigl[fX\extsubseteq X\,\then\,Y\extsubseteq X\bigr]
\end{array}$$ %}\end{Definition}
When $f$ satisfies $\Clop(f,A)$, we call $f$ a {\em classification operation
on $A$}.
When $f$
satisfies $\Clop(f,A)$ and $\Ext(f,A)$, we call $f$ {\em extensional} or an
{\em extensional
operation on $A$}. When $f$ satisfies $\Clop(f,A)$ and $\Mon(f,A)$, we say
that
$f$  is a {\em monotone operation on $A$}.
Since monotonicity entails extensionality, a monotone operation is always
extensional.

 Now we state %$\MID_A$ and
$\UMID_A$.\\[2ex]
{\bf $\MID_A$ (Monotone Inductive Definition on $A$)}
$$\forall f\,[\Clop(f,A)\land\Mon(f,A)\then\exists Y\,\Lfp({Y},{f,A})].$$
{\bf ${\mathbf{UMID}}_A$ (Uniform Monotone Inductive Definition on $A$)}
$$\forall f\,[\Clop(f,A)\land \Mon(f,A)\then\Lfp(\clfp(f),f,A)].$$
$\UMID_A$ states that if $f$ is monotone on subclassifications of $A$, then
$\clfp(f)$ is a least fixed point of $f$.

Let $\V{}$ be the universe, i.e. $\V{}:=\{x: x=x\}$.
By $\MID$ and $\UMID$ we denote
the principles $\MID_{\V{}}$ and $\UMID_{\V{}}$, respectively.\footnote{The
acronym for the principle $\MID$ in Feferman's paper \cite{Fef79}, section 7
was $\mbox{MIG}\upharpoonright$.}  }\end{deff}

The strength of the various classical versions
was determined as a result of several papers \cite{Ra96,Ra98,Ra99,grs}.
The $\MID$ case is dealt with in \cite{grs,Ra02}. \cite{Ra02} provides a survey of all known results in the classical case.
 $\UMID_{\mathbb{N}}$ was shown to be related to
subsystems of second order arithmetic based on
${\Pi}^1_2$ comprehension.

To relate the state of the art in these matters we shall need some
terminology. Below  we shall distinguish between the classical and
the intuitionistic version of a theory by appending the
superscript $c$ and $i$, respectively. For a system $S$ of
explicit mathematics we denote by $S\!\restriction$ the version
wherein the induction principles for the natural numbers and for
inductive generation are restricted to sets. $\indd$ stands for
the schema of induction on natural numbers for arbitrary formulas
of the language of explicit mathematics. $\pizweir$ denotes the
subsystem of second order arithmetic (based on classical logic)
with ${\mathbf{\Pi}}^1_2$-comprehension but with induction
restricted to sets, whereas $\pizwei$ also contains the full
schema of induction on $\bN$.

 \cite{Ra98,Ra99} yielded the following results:
\begin{thm}\label{Michael}
\begin{itemize}
\item[(i)] $\pizweir$ and $\tnurc$ have the same proof-theoretic
strength. \item[(ii)] $\pizwei$ and $\tnurindc$
 have the same proof-theoretic strength.
\end{itemize}
\end{thm}
The first result about $\UMID_{\bN}$ on the basis of
intuitionistic explicit mathematics was obtained by Tupailo in \cite{Tumid}.
\begin{thm}\label{Sergei}
$\pizweir$ and $\tnuri$ have the same proof-theoretic strength.
\end{thm}
\cite{Tumid} uses a characterization of $\pizweir$ via a classical
$\mu$-calculus (a theory which extends the concept of an inductive
definition), dubbed ${\mathbf{ACA}}_0({\mathcal L}^{\mu})$, given
by M\"ollerfeld \cite{Moe02} and then proceeds to show that
${\mathbf{ACA}}_0({\mathcal L}^{\mu})$ can be interpreted in its
intuitionistic version, ${\mathbf{ACA}}_0^i({\mathcal L}^{\mu})$,
by means of a double negation translation. Finally, as the latter
theory is readily interpretable in $\tnuri$, the proof-theoretic
equivalence stated in Theorem \ref{Sergei} follows in view of
Theorem \ref{Michael}.

The proof of \cite{Tumid}, however, does not generalize to
$\tnurindi$ and extensions by further induction principles. The
main reason for this is that adding induction principles such as
induction on natural numbers for all formulas
 to  ${\mathbf{ACA}}_0({\mathcal L}^{\mu})$
only slightly increases the strength of the theory and by no means
reaches the strength of $\pizwei$. In order to arrive at a
$\mu$-calculus of the strength of $\pizwei$ one would have to
allow for transfinite nestings of the $\mu$-operator of length
$\alpha$ for any ordinal $\alpha<\varepsilon_0$. As it seems to be
already a considerable task to get a clean syntactic
formalization of transfinite $\mu$-calculi (let alone  furnishing
double negation translation thereof), this paper will proceed
along a different path. In actuality, much of the work was already
accomplished in \cite{Ra98}, where it was shown that
 $\pizweir$ and  $\pizwei$ can be reduced to operator theories
 $\topp$ and $\Topp$, respectively.
A careful axiomatization of the foregoing theories in conjunction
with results from \cite{Ra96} showed that they lend themselves
to double negation translations and thus can be translated  into
their intuitionistic counterparts. As the intuitionistic theories
can be easily viewed as subtheories of $\tnuri$ and
$\tnurindi$, respectively, one can conclude the
following result.
\begin{theorem}\label{MT}
\begin{itemize}
\item[(i)] $\pizweir$ and $\tnuri$ have the same proof-theoretic
strength. \item[(ii)] $\pizwei$ and $\tnurindi$
 have the same proof-theoretic strength.
\end{itemize}
\end{theorem}
\prf See \cite{rtu}. \qed
Through Theorem \ref{MT} one also gets a different proof of Theorem
\ref{Sergei} which does not hinge upon \cite{Moe02}.

\begin{rmk}{\em Virtually nothing is currently known about the strength of $\TO^i+\MID$ and variants. In the classical case
there is a close relationship with parameter-free $\Pi^1_2$-comprehension.
It would be very interesting to investigate whether the strength of $\MID$ diminishes in the intuitionistic setting.}
\end{rmk}

The strength of explicit mathematics with principle like $\UMID_{\mathbb N}$ and even $\MID$ considerably exceeds that of Martin-L\"of type theory.
This has a bearing on foundational questions such as the limit of constructivity or the limits of different concepts
of constructivity.
 In \cite{rsynthese,limits} an attempt is made to delineate the form of constructivism underlying Martin-L\"of type theory,
 suggesting that $\TO^i+\UMID_{\mathbb N}$ lies beyond its scope.

\paragraph{Acknowledgement}
Part of the material is based upon research supported by the EPSRC of the UK through grant No. EP/K023128/1. This research was also supported by a Leverhulme Research Fellowship and a Marie Curie International Research Staff Exchange
 Scheme Fellowship within the 7th European Community Framework Programme.
 This publication was also made possible through the support of a
grant from the John Templeton Foundation.

 Thanks are also due for the invitation to speak at the American Annual Meeting of the Association for Symbolic Logic (University of Connecticut, Storrs, 23 May, 2016) where the material of the first seven sections was presented.

\end{document}